\documentclass[11pt]{extarticle}
\usepackage[format=plain, font=footnotesize]{caption}
\usepackage{amsmath, amsthm, amssymb, mathtools}
\usepackage[hypertexnames=false,pdftex,backref=page,
 	pdfpagemode=UseNone,
 	breaklinks=true,
 	extension=pdf,
 	colorlinks=true,
 	linkcolor=blue,
 	citecolor=red,
 	urlcolor=blue,
 ]{hyperref}

\usepackage[shortlabels]{enumitem}
\usepackage{graphicx}
\usepackage{tikz}
\usepackage{xcolor}
\usepackage{url}
\usepackage{multirow}
\usepackage{chemfig} 
\usepackage[version=3]{mhchem}
\usepackage{cleveref}
\usepackage{mlalgebra}
\hypersetup{linktocpage,linktoc=all,
  colorlinks=true,
  citecolor=blue,
  linkcolor=blue,
  urlcolor=black
}
\usepackage{tocloft}

\usepackage[numbers,sort]{natbib}

\newcommand{\RR}{\mathbb{R}}
\newcommand{\CC}{\mathbb{C}}
\newcommand{\PP}{\mathbb{P}}

\DeclareMathOperator{\tr}{tr}

\DeclareMathOperator{\sign}{sign}
\DeclareMathOperator{\rank}{rank}

\definecolor{darkgreen}{rgb}{0.0, 0.3, 0.0}
\definecolor{darkblue}{rgb}{0.0, 0, 0.9}

\newcommand{\xx}{\mathbf{x}}
\newcommand{\yy}{\mathbf{y}}
\oddsidemargin \evensidemargin
\let\oldbibliography\thebibliography
\renewcommand{\thebibliography}[1]{%
  \oldbibliography{#1}%
  \setlength{\itemsep}{0pt}%
}

\newtheorem{theorem}{Theorem}
\numberwithin{theorem}{section}
\newtheorem{proposition}[theorem]{Proposition}
\newtheorem{lemma}[theorem]{Lemma}

\theoremstyle{definition}
\newtheorem{definition}[theorem]{Definition}

\newtheorem{question}[theorem]{Question}

\crefname{equation}{}{}
\crefname{equation}{}{}
\crefname{figure}{Figure}{Figure}
\crefname{section}{Section}{Section}
\crefname{lemma}{Lemma}{Lemma}
\crefname{proposition}{Proposition}{Proposition}
\crefname{theorem}{Theorem}{Theorem}
\crefname{corollary}{Corollary}{Corollary}
\crefname{definition}{Definition}{Definition}
\crefname{notation}{Notations}{Notation}
\crefname{remark}{Remark}{Remark}
\crefname{claim}{Claim}{Claim}
\crefname{assumption}{Assumption}{Assumption}

\title{Nonlinear Algebra and Applications}
\date{}
\author{
Paul Breiding, T\"urk\"u \"Ozl\"{u}m \c{C}elik,
Timothy Duff,
Alexander Heaton,
Aida Maraj,\\
Anna-Laura Sattelberger,
Lorenzo Venturello, and O\u{g}uzhan Y\"ur\"uk
}

\begin{document}
\maketitle

\begin{abstract}
We showcase applications of nonlinear algebra in the sciences
and engineering. Our review is organized into eight themes: polynomial
optimization, partial differential equations, algebraic statistics,
integrable systems, configuration spaces of frameworks, biochemical reaction
networks, algebraic vision, and tensor decompositions. Conversely, developments
on these topics inspire new questions and algorithms for algebraic geometry.
\end{abstract}

\section*{Introduction}
This is a review about nonlinear algebra and applications. We are eight researchers working on different aspects of nonlinear algebra. In this article, we aim to give a glimpse of the diversity of the field. The word ``applications'' refers to the fact that each of us connects their research to other parts of mathematics or concrete problems arising in the sciences.

So what {is} nonlinear algebra? In simple terms, nonlinear algebra adds the prefix ``non'' to linear algebra. But it is not just an extension that moves the focus from linear to nonlinear equations. Nonlinear algebra encompasses many different methods from fields like combinatorics, group theory, topology, convex and discrete geometry, tensors, number theory, representation theory, and algebraic geometry. All of this is combined with computational methods using symbolic and numerical computations to tackle concrete problems. We quote the following paragraph from the preface of the book \emph{Invitation to Nonlinear Algebra}~\cite{MS2020} by Micha\l{}ek and Sturmfels:

\begin{quote}
\emph{``Linear algebra is the foundation of much of mathematics, particularly in applied mathematics. Numerical linear algebra is the basis of scientific computing, and its importance for the sciences and engineering can hardly be overestimated. The ubiquity of linear algebra masks the fairly recent growth of nonlinear models across the mathematical sciences. There has been a proliferation of methods based on systems of multivariate polynomial equations and inequalities. This is fueled by recent theoretical advances, efficient software, and an increased awareness of these tools. Application areas include optimization, statistics, complexity theory, among many others.''}
\end{quote}
The lecture~\cite{MS2018} based on~\cite{MS2020} is meant to provide an introduction to the basic theory underlying the many different fields within nonlinear algebra. Our article can be seen as supplementary material showing examples of how this theory is applied. With our article we invite researchers from the applied mathematical sciences to familiarize themselves with the ideas and concepts of our field and convince them that nonlinear algebra provides powerful tools and methods that can help solve their problems. While linear algebra is well-established in applications and broadly used because of its accessibility, {\em non}-linear algebra might,  to people that are unfamiliar with algebraic geometry, appear intimidating and possibly too abstract for being applicable.
We would like to borrow the following words by Dolotin and Morozov~\cite{DM2007}: {\em ``Non-linear algebra is as good and as powerful as the linear one.''} In the sequel of this article, we hope to dissipate fear of contact and convince you that nonlinear algebra complements and enlarges the toolbox of linear algebra. 

\bigskip

\noindent \textbf{Outline.} After a short section recalling basic definitions from algebraic geometry, the article is subdivided into eight different sections---one for each author. Each of us explains their view on nonlinear algebra and how it relates to other fields and applications. Yet, we  emphasize that our overview is not comprehensive. It is not our intention to map out all of nonlinear algebra. We  seek to give a short glimpse of the diversity of nonlinear algebra and convince you of its elegance and utility in concrete applications.

\begingroup\hypersetup{linkcolor=black}\normalfont
\tableofcontents\endgroup

\setcounter{section}{-1}
\section{Some basics about algebraic varieties}\label{section intro}
In nonlinear algebra, models are defined by polynomials.
For instance, in many applications our model~$V$ may be a real algebraic variety, defined by polynomial equations, i.e.,
\begin{equation}\label{def_variety}
V = \{\mathbf{x}\in\mathbb R^n \mid f_1(\mathbf{x}) = \cdots = f_s(\mathbf{x})=0\},
\end{equation}
where $f_1,\ldots,f_s$ are polynomials in $n$ variables $x_1,\ldots,x_n$.
In many branches of nonlinear algebra, one moves from the real numbers to the complex numbers. A fundamental motivation for this is that the field of complex numbers is algebraically closed.
One then considers \emph{complex} algebraic varieties
\begin{equation}\label{def_variety_complex}
V^{\mathbb C}= \{\mathbf{x}\in\mathbb C^n \mid f_1(\mathbf{x}) = \cdots = f_s(\mathbf{x})=0\}.
\end{equation}
In fact, one can go even one step further and consider complex \emph{projective} varieties. Those are algebraic varieties in the complex projective space~$\mathbb P^n = \{\mathbb C \cdot \mathbf{y}\mid \mathbf{y}\in\mathbb C^{n+1}\setminus \{{\bf 0}\}\}$, which is defined as the set of complex lines through the origin in $\mathbb C^{n+1}$. Those varieties are cut out by a family of homogeneous polynomials in the variables $x_0,\ldots,x_n.$
Complex algebraic varieties of the form~\eqref{def_variety} embed into complex projective $n$-space via the embedding $
\iota \colon  \mathbb C^n \to \mathbb P^n,$ $\mathbf{x} \mapsto \mathbb C\cdot (1,\mathbf{x}),
$
i.e., $\mathbf{x}\in \mathbb C^{n}$ is mapped to the line through the point in $\mathbb C^{n+1}$ which is obtained by appending a~$1$~to the vector~$\mathbf{x}$.
Associated to a projective variety $X\subseteq \mathbb P^n$ we have its underlying variety in~$\CC^{n+1}$, which we denote by~$\widehat{X}$. The latter is a cone, since $X$ is cut out by homogeneous polynomials. The motivation for studying the projective space and its subvarieties is that $\mathbb P^n$ is compact.
While algebraic subvarieties of~$\mathbb{P}^n$ are called {\em projective}, varieties in~$\mathbb C^n$ are called \emph{affine} or \emph{quasi-projective}.

The set of complex algebraic varieties in $\mathbb C^n$ is closed under taking intersections and under taking finite unions. Thus, they form a system of closed sets in a topology on~$\mathbb C^n$. This topology is called the \emph{Zariski topology}. We can also define the Zariski topology on $\mathbb R^n$ or $\mathbb P^n$ using real or projective varieties instead.

Associated to a complex variety $V^\mathbb{C}$ is its ideal $I(V^\mathbb{C})$, which is defined as the ideal in the polynomial ring $\mathbb{C}[\mathbf{x}]=\mathbb C[x_1,\ldots,x_n]$ containing all polynomials that vanish on $V^\mathbb{C}$. Recall that an {\em ideal} $I$ of $\mathbb{C}[\mathbf{x}]$ is a nonempty subset of polynomials 
that is closed under addition and under multiplication by elements in $\mathbb{C}[\mathbf{x}]$, which is denoted by $I \triangleleft \mathbb C[\mathbf{x}]$. Studying geometric properties of a variety translates to studying ideals
from the perspective of commutative algebra. This approach constitutes the field of algebraic geometry  \cite{Hartshorne,Liu,Vakil}. For an introduction to commutative algebra, we refer our readers to the textbook~\cite{Eisenbud}. By Hilbert's Nullstellensatz~\cite[Chapter 6]{MS2020}, for any ideal $I$ of $\mathbb{C}[\mathbf{x}]$ there exists a finite set of polynomials that generates~$I$. This gives a representation of $I$ that is not unique.
Some of such representations of an ideal are better suited for computations than others. One representation that is particularly useful in applications is in terms of a Gr\"obner basis~\cite{whatisGrobner,CLO,SingularCommAlg,sturmfels1996grobner,Sturmfels02solvingsystems}. Many important invariants, such as the dimension or the degree of the variety, can be directly read from a Gr\"obner basis. Another important representation is the primary decomposition, which contains information about  the irreducible components of the variety, i.e., how the variety decomposes set-theoretically.

In order to define the notion of genericity---which is crucial for certain results in nonlinear algebra---let~$X$ be an irreducible complex or projective variety. We call a property \emph{generic} for~$X$, if it holds for all points of~$X$ outside of an algebraic subvariety $\Delta\subsetneq X$. The subvariety $\Delta$ is called the \emph{discriminant} for the property. We say that $p\in X$ is \emph{generic} if $p\not\in\Delta$. This terminology refers to the fact that, in case of a generic property, almost all points of the variety are generic. An applied scientist would rephrase ``almost all'' as ``with probability $1$''.
Let us illustrate the definition in a concrete example: a quadratic polynomial of the form $f = ax^2+bx+c\in \CC[x]$ has exactly two distinct zeros in~$\mathbb C$ if and only if $b^2 - 4ac\neq 0$. The discriminant $\Delta = \{ b^2 - 4ac = 0\}$ is an algebraic variety in the space of all quadratic polynomials. In this example, the definition of~$\Delta$ is clear. In other contexts, the definition of~$\Delta$ usually is less obvious and equations might not be known. Nevertheless, in the following sections we will use the term generic, even if we don't know~$\Delta$ explicitly, but know it exists. Moreover, when we call a point {\em generic} in the following, we implicitly mean generic with respect to a property given by the context. In several applications, one is interested in the question of whether or not a set of real points $R\subseteq X$ contains generic points. For this it often suffices to show that the Zariski closure of $R$ is~$X$.

Computations are an important cornerstone of nonlinear algebra. Researchers in this field use a collection of different methods ranging from numerical algorithms such as {polynomial homotopy continuation}~\cite{Sommese:Wampler:2005} to symbolic methods based on Gr\"obner bases  \cite{Sturmfels02solvingsystems} or  combine symbolic and numerical algorithms~\cite{HL2015}. For carrying out symbolic computations, various computer algebra systems are available, such as {\tt Macaulay2}~\cite{M2}, {\tt SageMath}~\cite{sagemath}, or {\tt Singular}~\cite{Singular}, to mention a small selection. Numerical methods are implemented in several sophisticated software packages like~\cite{Bertini,Hom4PSArticle,HC.jl,NumericalAlgebraicGeometryArticle,PHCpack}.


\section{Polynomial optimization}\label{section polopt}
\emph{by Lorenzo Venturello}
\medskip

In this section, we consider the problem of optimizing a polynomial function $f:\mathbb{R}^n\to \mathbb{R}$ on a subset of $\mathbb{R}^n$ defined by a finite collection of polynomial equations and inequalities. A set of this form is called a \emph{basic semialgebraic set}. We hence consider the following polynomial optimization problem:
\begin{align}
\begin{split}
	\label{eq:problem}   
	{\rm Minimize} \,\,\,
	f(x_1,\dots,x_n) \quad  
	\text{subject to}\,\,\, & f_i(x_1,\dots,x_n)=0 \,\,\,\,
	\text{for all} \,\, \,1\leq i\leq k,\\
	&g_j(x_1,\dots,x_n)\geq0 \,\,\,
	\text{for all} \,\, \,1\leq j\leq m, 
	\end{split}
	\end{align}
with $f_i,g_j\in \mathbb{R}[x_1,\dots,x_n].$ This very general setting specializes to a number of relevant convex optimization problems in mathematics, such as {linear programming}, {semidefinite programming}, and {conic optimization}. For the purpose of this survey, we will assume that $m=0,$ i.e., no inequalities are considered. The more general case can be essentially derived with a case study depending on which inequality constraints are active. 
	
Our goal is to employ methods from nonlinear algebra to quantify the complexity of computing an exact solution of \eqref{eq:problem}. This leads to the discovery of intriguing connections with classical topics in algebraic geometry. We then observe that if a smooth point $\mathbf{x}^*=(x_1^*,\dots,x_n^*)$ of the variety defined by the constraints~$f_1(\mathbf{x})=\dots = f_k(\mathbf{x})=0$ is an optimal solution to~\eqref{eq:problem}, then it must satisfy the \emph{Karush--Kuhn--Tucker} (KKT) conditions, which are given by
	\begin{align}\label{eq: KKT}
		\nabla f(\mathbf{x}^*) - \sum_{i=1}^{k} \lambda_i \nabla f_i(\mathbf{x}^*)=0 \quad \text{and} \quad
		f_j(\mathbf{x}^*)=0 \,\,\, \text{for all} \,\, \, 1\leq j \leq k,
	\end{align}
for some $\lambda_i\in\mathbb{R}.$ This implies that smooth optimal solutions to \eqref{eq:problem} are to be found among the solutions of the system of polynomial equations \eqref{eq: KKT}. From now on we move to the field of complex numbers and hence consider also complex solutions to \eqref{eq: KKT}. As it has been discussed in \Cref{section intro}, this is a key technical condition to employ algebraic methods, as it allows us to study the \emph{generic} instance of our problem. We now consider a fundamental invariant of a polynomial optimization problem.
\begin{definition}
	A point $\mathbf{x}^*\in\mathbb{C}^n$ such that $(\mathbf{x}^*,\mathbf{\lambda}^*)$ is a solution of \eqref{eq: KKT} for some $\mathbf{\lambda}^*\in\mathbb{C}^k$ is a \emph{critical point} of \eqref{eq:problem}. We define the \emph{algebraic degree} of an optimization problem to be the number of critical points counted with multiplicity, when it is finite.
\end{definition}
In general the critical points of \eqref{eq:problem} can form a positive-dimensional variety. For example, if the polynomials $f$ and $f_i$ are homogeneous, it is not hard to see that if $\mathbf{x}^*$ is a critical point, then $\mu\mathbf{x}^*$ is also a critical point for every $\mu\in\mathbb{C}.$ However, these are rather special instances and from now on we shall assume the existence of a finite number of critical points.

The algebraic degree can be described also in another way. The coordinates of smooth optimal solutions are \emph{algebraic functions} in the input, i.e., each coordinate is the zero of a univariate polynomial in the coefficients of $f$ and of the $f_i$'s. The degree of this polynomial is the algebraic degree of the corresponding optimization problem.
The next result provides a closed formula for the algebraic degree when $f_1,\dots,f_k$ are generic polynomials with~$k\leq n.$ Here the word ``generic'' indicates that there is an open Zariski dense set of polynomials satisfying this property. Furthermore, the techniques to obtain these results often rest on projective geometry.
\begin{theorem}[{\cite{NR}}]\label{thm: nr}
	Let $k\leq n$ and let $d=\deg(f)$ and $d_i=\deg(f_i)$ for every $i=1,\dots,k.$ If~$f_1,\dots,f_k$ are generic,  then the algebraic degree of \eqref{eq:problem} is equal to 
	\[
		d_1\cdots d_k\sum_{i_0+\dots+i_k=n-k} (d-1)^{i_0}(d_1-1)^{i_1}\cdots (d_k-1)^{i_k}.
	\]
	If $f_1,\dots,f_k$ are not generic, then the formula gives an upper bound for the algebraic degree of \eqref{eq:problem}.
\end{theorem}
For instance, if we minimize a linear function over the vanishing locus of a generic polynomial $f_1\in\mathbb{R}[x_1,\dots,x_n]$ of degree $d_1,$ the formula in \Cref{thm: nr} states that we should expect $d_1(d_1-1)^{n-1}$ critical points. Another interesting case is when $k=n.$ In this case (and when the polynomials are generic) the feasibility region of \eqref{eq:problem} is zero-dimensional and, by B\'{e}zout's theorem, it consists of $d_1\cdots d_k$ many points. This is precisely the value predicted by the formula above.

Similar techniques have been applied to study the algebraic degree of semidefinite programming \cite{NRS, vBR, MS2020}, the problem of optimizing the Euclidean distance over an algebraic variety \cite{DHOST, OtSo}, and maximum likelihood estimation in algebraic statistics \cite{algstat2018}, which will be further discussed in Section \ref{section AlgStat}. In those settings, the nature of the problem imposes a certain structure on the polynomials $f_i,$ far from the generic behavior. This asks for a refinement of \Cref{thm: nr}. Once again we can make use of tools known to algebraic geometers. For example, every smooth projective  variety $X\subseteq \mathbb{P}^n$ comes with a sequence of nonnegative numbers $(\delta_0(X),\dots,\delta_{n-1}(X))$ called \emph{polar degrees} of $X.$ The number $\delta_r(X)$ counts the number of intersection between the {conormal variety} of $X$ in $\mathbb{P}^n\times \mathbb{P}^n$ \cite[Chapter 1, Section 3]{GKZ} and a pair of generic projective linear spaces $(L,M),$ with $\dim L = n-r$ and $\dim M = r+1.$
\begin{theorem}[{\cite[Section 5]{DHOST}, \cite[Theorem 13]{CJMSV}}]\label{thm: polar} 
	Let $f, f_1,\dots,f_k\in\mathbb{R}[x_1,\dots,x_n]$. Let $r\geq 0$ be such that $1\leq n-r \leq k$ and assume that $f_1,\dots,f_{n-r}$ are generic linear polynomials. Let $X$ be the projective closure of the variety defined by $f_{n-r+1}(\mathbf{x})=\cdots= f_k(\mathbf{x})=0.$ Then the algebraic degree of problem \eqref{eq:problem} is bounded above by $\delta_r(X),$ the $r$-th polar degree of $X.$ 
\end{theorem}
In other words, the polar degrees of (the projective closure of) a variety are precisely the algebraic degree of optimizing a generic linear functional over the intersection of the variety with a generic linear space. The appropriate polar degree is determined by the dimension of this linear space.
Closed formulas for the polar degrees are known for a number of interesting varieties. In general, they can be obtained from explicit computations in the appropriate {Chow ring} or with the help of a computer algebra system such as {\tt Macaulay2}~\cite{M2}, as described in \cite[Section 5]{DHOST}.

For a simple example of how this result is finer than \Cref{thm: nr}, let $f_1\in\mathbb{R}[x_{12},x_{13},x_{21},x_{22},x_{23}]$ be a generic linear polynomial and let $f_2=x_{22}-x_{12}x_{21},$ $f_3=x_{23}-x_{13}x_{21}$ and $f_4=x_{12}x_{23}-x_{13}x_{22}.$ The three  polynomials of degree $2$ are the $2\times 2$-minors of a generic $2\times 3$-matrix with $x_{11}=1$ and they cut the variety of such $2\times 3$-matrices of rank at most $1.$ We want to optimize a generic linear function over the smooth variety of $\mathbb{R}^5$ cut out by $f_1,\dots,f_4.$ In the case of generic polynomials with $d=1,$ $d_1=1$ and $d_2=d_3=d_4=2,$ the formula in \Cref{thm: nr}  predicts an algebraic degree equal to $24.$ On the other hand, \Cref{thm: polar} states that this algebraic degree is bounded above by the first polar degree of the projective variety of  $2\times 3$-matrices of rank $1,$ which is equal to~$3.$ This difference is precisely due to the fact that $f_2,f_3,f_4$ are definitely not generic, and they come from a very structured problem. Moreover, \Cref{thm: polar} gives only an upper bound: using {\tt Macaulay2}~\cite{M2}, we can verify that for this example the system \eqref{eq: KKT} has only two solutions. This second discrepancy is due to the fact that solutions might lie on the so-called \emph{hyperplane at infinity}, a feature which occurs in the transition between affine and projective spaces \cite[Chapter 8]{CLO}. An even more peculiar behavior can be observed if we remove $f_1$ from this example and optimize a generic linear function on the variety defined by $f_2,f_3$, and $f_4.$ In this case, the system \eqref{eq: KKT} has no solution---not even over the complex numbers. This has a beautiful explanation via the theory of dual varieties and it corresponds to the fact that the dual of the variety of $2\times 3$-matrices of rank~$1$ is not a hypersurface~\cite{GKZ}.  

An intriguing application is to describe the algebraic degree of computing the distance from a point to an algebraic variety using a polyhedral norm. In this case, we have a finite number of optimization problems of the form of \Cref{thm: polar}; one for every $(r-1)$-dimensional face of the unit ball. As an example, consider the so-called \emph{Wasserstein distance} induced by a metric on a finite space, which is a main concept in optimal transport and machine learning. The algebraic treatment of this optimization problem has been described recently in \cite{CJMSV, macis}.

\section{Partial differential equations}\label{section PDEs}
\emph{by Anna-Laura Sattelberger}

\medskip

Algebraic analysis investigates linear partial differential equations (PDEs) with tools from algebraic geometry, category theory, complex analysis, and noncommutative Gr\"obner bases by encoding linear PDEs as $D$-modules~\cite{HTT08, Bj2, Coutinho}. The latter arise in  Hodge theory, microlocal calculus, mirror symmetry, optimization, representation theory, statistics, and many more fields of mathematics. In its simplest form---which is at the same time the most prominent one in applications---$D$ denotes the {\em Weyl algebra}. It is a noncommutative algebra over the complex numbers, gathering linear partial differential operators with polynomial coefficients. In formal terms, the Weyl algebra $$D \, \coloneqq \, \CC[x_1,\ldots,x_n]\langle \partial_1,\ldots,\partial_n \rangle$$ is the free algebra over the complex numbers generated by the variables $x_1,\ldots,x_n$ and partial derivatives $\partial_1,\ldots,\partial_n,$ modulo the following relations: all generators are assumed to commute except $x_i$ and $\partial_i.$ Their  commutator $[\partial_i,x_i]=\partial_ix_i-x_i\partial_i$ is~$1,$ not $0.$ This identity encodes Leibniz' rule in a formal way. Since the Weyl algebra is noncommutative, one needs to distinguish between left and right multiplication. One usually studies {\em left} $D$-ideals and $D$-modules, i.e., multiplication by elements of $D$ is from the left.  Systems of homogeneous partial differential equations are then encoded as $D$-ideals, or, more generally, as $D$-modules. 
In analogy to encoding an algebraic number as the root of a polynomial, a function which is the solution of a sufficiently large system of linear PDEs can be encoded by a $D$-ideal together with finitely many initial conditions. Those functions are called {\em holonomic} and were first studied by Zeilberger~\cite{Zei90}. They fulfill certain closure properties and can be studied in terms of Gr\"obner basis computations in the Weyl algebra \cite{SST00}. Many functions in the sciences are indeed holonomic, such as rational functions, some trigonometric functions, hypergeometric functions,  Airy's functions and some  more special functions, and many more. In~\cite{SatStu19}, the authors give an introduction to $D$-modules and holonomic functions. As demonstrated therein, $D$-modules have a strong computational flavor and are useful in many applications in the sciences and algebraic geometry---as well as the other way round. Tools from algebraic geometry can help to explicitly construct solutions to a system of differential equations encoded by a \mbox{$D$-ideal}; for {Frobenius} ideals, solutions can be constructed in terms of primary decompositions of the underlying~ideal (cf. \cite[Section~2.3]{SST00}). 

Computations around $D$-ideals and holonomic functions are often cumbersome to carry out by hand. Several software packages are available for carrying out those computations, for instance in {\tt Macaulay2}~\cite{M2}, {\tt Magma}~\cite{magma} (e.g. the package {\tt periods}~\cite{periods}), {\tt Maple}~\cite{Maple} (e.g. the package {\tt gfun}~\cite{gfun}), {\tt Mathematica} (e.g. the {\tt HolonomicFunctions}~\cite{CK} package), {\tt SageMath}~\cite{sagemath} (e.g. the package {\tt ore\_algebra}~\cite{orealgebra}), or {\tt Singular:Plural}~\cite{Singular}.

The {\em principal symbols} of the differential operators contained in a $D$-ideal, i.e., the part containing the derivatives of highest orders, cut out the {\em characteristic variety} of the $D$-ideal. This is an algebraic variety in $\mathbb{C}^{2n}$ from which one can read the so-called {\em singular locus} of the $D$-ideal; there, complex analytic solutions to the system of differential equations encoded by the $D$-ideal might have singularities. 

In a more general setup, one studies the category of sheaves of $D$-modules on a complex manifold or a smooth algebraic variety. This abstract language allows for deep structural insights to linear PDEs. For instance, the regular Riemann--Hilbert correspondence, proven in the 1980s by Kashiwara~\cite{Kas84} and Mebkhout~\cite{Meb84} independently, gives a refined answer to Hilbert's $21$st problem in a generalized setup. It enables to replace systems of linear PDEs with ``mild'' singularities by their topological counterpart derived from their holomorphic solutions and vice versa.

As already mentioned, the theory of $D$-modules comes with a vast range of applications. Among them are the high-precision computation of periods~\cite{periods} or the volume of compact semi-algebraic sets~\cite{LMS19}, and the {holonomic gradient method} (HGM). This method is a numerical scheme for the evaluation of holonomic functions~\cite{NNNOSTT11}.
Exploiting the knowledge of an annihilating $D$-ideal, the evaluation of the gradient in each iteration step is reduced to a matrix multiplication. The holonomic gradient {\em descent} is a minimization scheme building on the HGM. For some functions, those methods are encoded in the software package {\tt hgm}~\cite{hgmR} in {\tt R}~\cite{R}. The said methods were  applied for instance to the inference of rotation data sampled according to the Fisher distribution \cite{SSTOT13, Koyama, ALSS} and to data arising from medical imaging. 

There are many more functions arising from statistics that suit this holonomic setup well. Muirhead \cite{Mui70,Mui82} observed that the cumulative distribution function of a Wishart matrix is related to the hypergeometric function ${_1F_1}$ of a matrix argument for which he constructed an annihilating $D$-ideal. This ideal was further studied in~\cite{HNTT13,Noro,GLS}. While hypergeometric functions of one complex variable are well-studied in terms of GKZ systems~\cite{surveyGKZ}, there is no intrinsic description of hypergeometric functions of a matrix argument in terms of~\mbox{$D$-modules} yet. For broadening the range of applications of $D$-modules in the sciences, one is in need of seeking out further functions that can be described by means of holonomic functions.

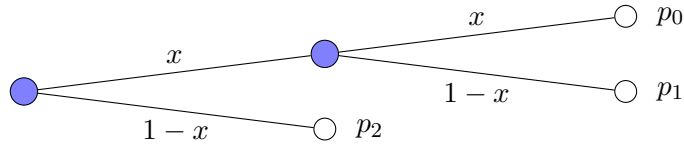
\begin{figure}
	\begin{center}
\begin{tikzpicture}[scale=1.0]
\node[fill,circle,draw=black, fill=blue!50, scale=1pt] (A) at (0,0) {};
\node[fill,circle,draw=black, fill=blue!50, scale=1pt] (B) at (4,0.5) {};
\node[fill,circle,draw=black, fill=white,scale=0.8pt] (C) at (4,-0.5) {};
\node[fill,circle,draw=black, fill=white,scale=0.8pt] (D) at (8,1) {};
\node[fill,circle,draw=black, fill=white,scale=0.8pt] (E) at (8,0) {};
\path [-] (A) edge node[above] {$x$} (B);
\path [-](A) edge node[below] {$1-x$} (C);
\path [-](B) edge node[above] {$x$} (D);
\path [-](B) edge node[below] {$1-x$} (E);
\node[right of=D] (p0) {\hspace*{-8mm}$p_0$};
\node[right of=E] (p1) {\hspace*{-8mm}$p_1$};
\node[right of=C] (p3) {\hspace*{-8mm}$p_2$};
\end{tikzpicture}
\caption{Staged tree~\cite{DMS21} modeling the discrete statistical experiment of flipping a biased coin twice.}
\label{fig:stagedtree}
\end{center}
\end{figure}

Via Bernstein--Sato ideals, $D$-modules are useful for the maximum likelihood estimation (MLE) problem in statistics. The {\em Bernstein--Sato polynomial} of a polynomial $f\in\CC[x_1,\ldots,x_n]$ is the monic polynomial~$b_f$ in $\CC[s]$ of smallest degree, s.t. there exists $P\in D[s]$ for which~\mbox{$ P\bullet f^{s+1} = b_f \cdot f^s.$} Here, $f^s$ denotes the symbolic power of~$f$ to the new variable~$s.$ The roots of $b_f$ are negative rational numbers~\cite{Kashiwara}. For a family of $k$ polynomials in $n$ variables, one studies the Bernstein--Sato {\em ideal}~\cite{SabbahBS} in the polynomial ring in $k$ variables. The original motivation of studying Bernstein--Sato polynomials was to construct a meromorphic continuation of the distribution-valued function~$f^s.$ Nowadays, Bernstein--Sato ideals are a prominent object of study in singularity theory~\cite{bvwz,Budur, BMT}.  For a concrete example, consider the following statistical experiment: {\em Flip a biased coin. If it shows head, flip again.} \Cref{fig:stagedtree} depicts the staged tree modelling this discrete experiment.

Let $(u_0,u_1,u_2)$ encode how many times each of the three states occurs when repeating the experiment several times. The MLE problem is to find the parameters that best explain this observed outcome; further details about maximum likelihood estimation are discussed in \Cref{section AlgStat}. By encoding statistical models as smooth subvarieties of projective space, likelihood geometry~\cite{HS14} provides algebraic tools for the computation around the MLE problem. 

\noindent For the outcome $(u_0,u_1,u_2),$ the MLE of the said experiment is
\begin{align}\label{MLEcoin}
\Psi(u_0,u_1,u_2) \,=\, \left( \frac{\left( 2u_0+u_1\right)^2}{\left( 2u_0+2u_1+u_2 \right)^2}, \,  \frac{\left( 2u_0+u_1\right) (u_1+u_2)}{\left( 2u_0+2u_1+u_2 \right)^2}, \, \frac{u_1+u_2}{2u_0+2u_1+2}  \right).
\end{align}
The triple $(x^2,x(1-x),1-x)$  is a parametrization of our statistical model. Its Bernstein--Sato ideal can be computed with the {\tt Singular} library {\tt dmod\_lib}~\cite{dmodlib}. 
As observed in~\cite{SatStu19}, the linear factors in the generator of the Bernstein--Sato ideal recover the occurring linear factors in the numerators  of the MLE in~\eqref{MLEcoin}. This observation is explained in a rigorous algebro-geometric way in~\cite{TropicalMLE} using tools from tropical geometry~\cite{MS15,OnTropComp} and Bernstein--Sato theory~\cite{Mai16}. The methods developed therein provide new tools for algebraic statistics and possibly particle physics; a link between MLE and scattering amplitudes was recently established by Sturmfels and Telen~\cite{ST20}.

Solutions to linear PDEs with {\em constant} coefficients can be constructed by the Ehrenpreis--Palamodov fundamental principle~\cite{Ehrenpreis,Palamodov}. In \cite{CRS21}, Cid-Ruiz and Sturmfels introduce {differential primary decompositions} for ideals in commutative rings, generalizing Noetherian operators in the setup of Ehrenpreis--Palamodov. On the repository website~\cite{mathrepo} of the Max Planck Institute for Mathematics in the Sciences (MPI MiS), they provide {\tt Macaulay2} code supplementing their work.

The study of PDEs by algebraic methods is not limited to the linear case. {\em Polynomial} differential equations are treated by differential algebra \cite{Ritt, Kolchin}. 
Among others, this setup can be used to study jet schemes of affine schemes~\cite{BMS}. In~\cite{AEMS21}, the authors transfer combinatorial aspects of commutative algebra~\cite{combalg} to differential algebra and link certain regular triangulations~\cite{triang} to a higher-dimensional analog of jets of the fat point~$x^p,$ building on the theory of differential Gr\"obner bases~\cite{Zobnin,Ollivier}. This case study suggests to further develop {combinatorial differential~algebra}. 

\section{Algebraic statistics}\label{section AlgStat}
\emph{by Aida Maraj}

\medskip

Algebraic statistics uses tools from algebraic geometry, commutative  algebra, combinatorics, and computational methods to solve  problems in probability theory and statistics. This relatively new area, marking its start with papers \cite{Pistone1996GeneralisedCW} in 1996 and  \cite{diaconis1998algebraic} in 1998,  provides new techniques for statistical problems, and has led to many interesting developments in mathematics. Sullivant's book \cite{algstat2018} contains a good sampling of algebraic statistics topics and their development. 

Given some independent, identically distributed observed data, the interest is in discovering features of the unknown probability distribution from which the data is drawn. \emph{Statistical models} are
collections of probability distributions possible for the data, typically parametrized or satisfying some property. Let  $\mathbf{u}\in \mathbb{Z}_{\geq 0}^n$ be the vector of counts for the given data with sample size $N=u_0+\cdots+u_{n-1}$. The vector $\mathbf{u}$ is often used to imply data with that vector of counts.
Given $\mathbf{u}$ and a statistical model $\mathcal{M},$ there are two questions to be asked: 
\begin{question}
\label{question1}
 How well does the model  $\mathcal{M}$ fit the data $\mathbf{u},$ i.e., what is the $p$-value for the  hypothesis that the unknown probability distribution for the given data is in $\mathcal{M}$?
\end{question}
\begin{question}
\label{question2}
 How to fit the model $\mathcal{M}$ to the data $\mathbf{u},$ i.e., what is the probability distribution in  $\mathcal{M}$ that best describes the data?
\end{question}
This section will discuss how Questions \ref{question1} and \ref{question2} are naturally translated to problems in algebraic geometry by identifying probability distributions with points in $\mathbb{R}^n$ and statistical models with (semi)algebraic sets. The focus will be on {discrete statistical models} and {multivariate Gaussian models}, for which algebraic methods have shown to be most useful. 

Discrete statistical models arise when random variables have a state space $\Omega$ of finite size, say~$n.$ This means that the experiments have $n$ possible outcomes. Assuming that each experiment has nonzero probability, discrete statistical models are subsets of the open
probability simplex 
\[\Delta_{n-1} \coloneqq \Big\{\mathbf{p}\in  \mathbb{R}_{>0}^{n}~|~ \sum\limits_{i=0}^{n-1} p_{i}=1 \Big\}.\] 
For instance, the staged tree model in \Cref{fig:stagedtree} is a discrete statistical model as it has a state space of size three. In this example, denote $p_0=P(X=\text{heads},Y=\text{heads}),$ $p_1=P(X=\text{heads},Y=\text{tails})$ and $p_2=P(X=\text{tails}),$ the probabilities for elements in the state space. The algebraic object for the staged tree model  are all points  $(p_0,p_1,p_2)\in \mathbb{R}_{>0}^3$ with $p_0+p_1+p_2=1$  that fulfill the equation $(p_0+p_1)p_1-p_0p_2=0.$ Hence our model is a curve in the triangle $\Delta_2.$ The ideal generated by the quadratic polynomial above is the \emph{vanishing ideal} of the model and  it encodes the staged tree model; see~\cite{duarte2020equations}. 

Algebra is particularly useful  if the vanishing ideal of a discrete statistical model is {toric}; see \Cref{thm:diastru}.
The ideal $I(\mathcal{M})$ for the discrete model $\mathcal{M}$ is {\em toric} if it is of form 
\(I(\mathcal{M})=\langle\mathbf{p}^{\mathbf{u}}-\mathbf{p}^{\mathbf{v}} ~|~ \mathbf{u},\mathbf{v}\in \mathbb{Z}_{\geq 0}^{n}, \mathbf{u}-\mathbf{v}\in \mathrm{ker}(A) \rangle,\) for some matrix $A\in \mathbb{Z}^{d\times n}$ (cf. \cite[Chapter~4]{sturmfels1996grobner}), where $\mathbf{p}^{\mathbf{u}}$ denotes $p_0^{u_0}\cdots p_{n-1}^{u_{n-1}}.$ Models with a toric vanishing ideal are called \emph{toric} or \emph{log-linear} and the design matrix $A$ is a \emph{sufficient statistic} for the toric model. In this case, the model itself can be restated as \[\mathcal{M} =\{\mathbf{p}\in \Delta_{n-1}~|~ (\mathrm{log} (p_0),\ldots,\mathrm{log} (p_{n-1}))\in \mathrm{rowspan}(A)\}.\]
{Markov bases} \cite{petrovic2019markov} of toric models  are critical in hypothesis testing; the Metropolis--Hasting algorithm uses them to produce Markov chains \cite[Section~9]{algstat2018}. Markov chains
are fundamental for computing the $p$-value for the null hypothesis that data $\mathbf{u}$ is drawn from a probability distribution in the given toric model \cite{aoki2012markov}, and hence answering \Cref{question1}.
\begin{theorem}[{\cite[Theorem~3.1]{diaconis1998algebraic}}]
\label{thm:diastru}
A collection of moves $\mathbf{b}_1,\dots,\mathbf{b}_m$ in $\mathbb{Z}^n_{\geq 0}$  is a Markov basis for the toric model $\mathcal{M}$ if and only if the set  $ \mathbf{p}^{\mathbf{b}^+_i}-\mathbf{p}^{\mathbf{b}^-_i},  \text{ for } 1\leq i\leq m,$
generates the toric ideal~$I(\mathcal{M}).$ 
\end{theorem}
A sufficient statistic  of toric models is often achieved effortlessly  from the monomial parametrization of the model in statistics; see for instance hierarchical models \cite{diaconis1998algebraic,hocsten2002grobner,maraj2020algebraic}, graphical models~\cite{geiger2006toric}, and balanced staged tree models \cite{duarte2020equations}. In other cases, a linear change of coordinates is needed to reveal the toric structure; see for instance conditional independence models and models arising from Bayesian networks \cite{garcia2005algebraic}, and some group based phylogenetics models \cite{sturmfels2005toric,coons2021toric}. The vanishing ideal for the staged tree model \Cref{fig:stagedtree} is toric in the new variables $q_0=p_0, q_1=p_0+p_1$ and $q_2=p_0+p_1+p_2$; see \cite{gorgen2021staged}.
The existence of matrix $A$ leads to involvement of polyhedral geometry \cite{bernstein2017unimodular,coons2020generalized}.

Likelihood estimation is the most common method for \Cref{question2}.  The  probability of observing data~$\mathbf{u}$ given the probability distribution $\mathbf{p}$ (see \cite[Section~1.1]{algstat2018})
is
\begin{align}\label{eq:discretelikeq}
   L(\mathbf{p}|\mathbf{u}) \coloneqq {\dfrac{N!}{\prod\limits_{i=0}^{n-1}(u_{i}!)}}\prod\limits_{i=0}^{n-1} p_{i}^{u_{i}}. 
\end{align}
The \emph{maximum likelihood estimate} (MLE) for data $\mathbf{u}$ in  the statistical model $\mathcal{M}$ is the probability distribution  $\hat{\mathbf{p}}$ in $\mathcal{M}$  so that under the assumed statistical model $\mathcal{M}$ the observed data is most probable. It is found by maximizing the likelihood function \ref{eq:discretelikeq} on points in $\mathcal{M}$. 
 A model $\mathcal{M}$ has \emph{rational maximum likelihood estimator} if for generic choices of $\mathbf{u},$ the MLE for $\mathbf{u}$ in $\mathcal{M}$ can be written as a rational function  in the entries of $\mathbf{u}$; see Example \Cref{MLEcoin} in \Cref{section PDEs}.
The multinomial coefficient and applying natural logarithm do not affect which probability is the maximizer. Hence, we optimize the \emph{log-likelihood function} $$\ell(\mathbf{p}|\mathbf{u})\coloneqq\sum\limits_{i=0}^{n-1} u_{i}\mathrm{log}(p_{i})$$ instead.
This is a linear function in terms of $\mathrm{log}(p_{i})$; compare this with optimizing a linear function in a variety in \Cref{section polopt}.
The MLE is found among critical points of $\ell_{\mathbf{u}}(\mathbf{p}).$ The \emph{maximum likelihood degree} (ML degree) of $\mathcal{M}$ is the number of complex critical points of $\ell_{\mathbf{u}}(\mathbf{p}),$ counted with multiplicity, for generic data $\mathbf{u}.$ Hence, the maximum likelihood degree measures the complexity of finding the MLE. 

For toric models, the ML degree is equal to  the number of intersection points of the toric model~$\mathcal{M}$ with the \emph{special} linear equations  $A(N\mathbf{p}-\mathbf{u})=0$ \cite[Corollary~7.3.9]{algcurves}, where $N$ is the number of observations.  Hence it is bounded above by the algebraic degree of the model. 
 Explicit computations can be done in {\tt Macaulay2} \cite{M2} via the package {\tt AlgebraicOptimization.m2}~\cite{harkonen2020algebraic}.

A toric model has rational MLE if and only if its ML degree is one \cite{huh2013varieties}.
Apart from special classes \cite{coons2021quasi,DMS21},  a classification of statistical models  with rational MLE with practical uses is to be found. Ultimately,  revealing the toric structure of a discrete statistical model, when it exists, is beneficial for answering Questions \ref{question1} and \ref{question2}.  


A multivariate Gaussian distribution $\mathcal{N}_n(\boldsymbol{\mu},\Sigma)$ is a continuous distribution determined by the mean vector $\boldsymbol{\mu}\in \mathbb{R}^n$ and the covariance matrix $\Sigma,$ which is a symmetric positive definite matrix. Denote by $\mathbb{S}_{>0}^n$ the space of symmetric $n\times n$ positive definite matrices. The parameters $\boldsymbol{\mu}$ and $\Sigma$ are independent of  each other. In computations, $\boldsymbol{\mu}$ is most commonly zero or the mean vector of the sample data, and $\Sigma$ is the challenging parameter to be found. Hence, in algebraic statistics,  multivariate Gaussian models are identified with sets of matrices in  $\mathbb{S}_{>0}^n.$  Such models include graphical models and colored graphical models \cite{SU2010}, Brownian motion tree models  \cite{Sturmfels2020BrownianMT}, and conditional independence models \cite{sullivant2009gaussian}. The inverse matrix of $\Sigma$ is the {\em concentration matrix} for the distribution. Given that both $\Sigma$ and the  concentration matrix $\Sigma^{-1}$ contain critical information about the statistics of the model, a main problem is to describe the model in both the space of covariance matrices, and its inverse space. The question has found partial answers when one of the sets is a linear system of symmetric matrices \cite{SU2010,Sturmfels2020BrownianMT}. 

The sample covariance matrix for  data  $\mathbf{u}_1,\dots,\mathbf{u}_d\in \mathbb{Z}_{\geq 0}^n$ is \[S=\dfrac{1}{d}\sum_{i=1}^d \mathbf{u}_i\mathbf{u}_i^T.\]  
The \emph{log-likelihood function} of a Gaussian model with  empirical covariance matrix $S$ is the linear function
$\ell_S: \mathbb{S}_{>0}^n \to~\mathbb{R}$ defined by $ \ell_S(M) \coloneqq \mathrm{log}\: \mathrm{det}(M) - \mathrm{trace}(S M).$ 
Maximizing $\ell_S$ on the space of concentration  matrices  produces the  {maximum likelihood estimator} and the {maximum likelihood degree} of the model. 
Maximizing $\ell_S$ on the set of covariance matrices gives the  \emph{reciprocal maximum likelihood estimator} and the \emph{reciprocal maximum likelihood degree} of the model. 

In the case of a linear concentration model $\mathcal{L},$ the ML degree of $\mathcal{L}$ is the number of covariance matrices $\Sigma$ in  $\mathcal{L}^{-1}$ satisfying linear equations arising from $\Sigma-S\in \mathcal{L}^{\perp}$ \cite[Corollary~7.3.10]{algstat2018}. This makes the algebraic degree of $\mathcal{L}^{-1}$ its upper bound. The references \cite{lauritzen2019maximum,Sturmfels02solvingsystems,SU2010} highlight some of the progress in this direction.
The packages {\tt SemidefiniteProgramming.m2}~\cite{cifuentes2020sums} and {\tt GraphicalModelsMLE.2} \cite{amendola2020computing} in {\tt Macaulay2} \cite{M2} are helpful for computations. 
Apart from special cases~\cite{dye2020maximum,fevola2020pencils}, a relation between the MLE and the reciprocal MLE is to be found.   

When one of the defining spaces of a Gaussian model is toric \cite{Sturmfels2020BrownianMT,misra2020gaussian},   one hopes to connect the model to  the associated toric discrete model. For instance, Brownian motion tree models are toric in the space of concentration matrices \cite{Sturmfels2020BrownianMT}. This same space defines a discrete toric model from phylogenetics. The reciprocal ML degree of the first model is  the ML degree of second model~\cite{boege2020reciprocal}. It is unexplored how tight the connections between the toric Gaussian model and the associated toric discrete model are and how to make use of these connections.  

Invariant theory has shown useful in computing maximum likelihood estimates and degrees~\cite{amendola2021invariant,derksen2020maximum,draisma2013groups}.  Another approach for tackling \Cref{question2} is to compute the point in a model that minimizes some distance of the sample point from the model \cite{CJMSV,DHOST}.  Model selection \cite[Chapter~17]{algstat2018} and identifiability \cite{allman2010identifiability,hollering2021identifiability} are other topics of current importance in algebraic statistics. See also \Cref{section tensors} for a related discussion about identifiability.

Lastly, the need of algebraic statistics to advance statistical questions has led to invention and development of tools in algebra; Huh's work \cite{huh2013varieties} on rational maximum likelihood estimates produced an interesting classification of varieties, and work of Hillar and Sullivant \cite{hillar2012finite} on finitely generated models up to a symmetric group action spawned the non-Noetherian theory in commutative algebra \cite{draisma2015noetherianity,krone2017hilbert,maraj2019equivariant,nagel2017equivariant,nagpal2021symmetric}. 


\section{Integrable systems}\label{section KP} 
\emph{by T\"urk\"u \"Ozl\"um \c{C}elik}
\medskip

Integrability of a system of partial differential equations (PDEs) reveals itself through the following features of the PDE: the essence of algebraic geometry, the presence of conserved quantities, and the existence of explicit solutions~\cite{HitSegWar}. It turns out that many known integrable systems can be obtained as reductions of the \emph{Kadomtsev--Petviashvili (KP) hierarchy}, which is an infinite set of nonlinear partial differential equations in infinitely many variables. This universal hierarchy consists of the {KP equation} and its infinitely many symmetries, for which the KP equation happens to be the initial member~\cite{Kod, NovManPitaZak}. Among the approaches to study the integrable systems is via bilinearizing the partial differential equation~\cite{Hir}. After this bilinearization, the integrable system is transformed to the so-called \emph{Hirota bilinear equation}~\eqref{eq:hirota}, which will be the principal object underlying the discussion about the KP equation in the sequel. 

The KP equation is a nonlinear partial differential equation that describes water waves. Its equation is written as follows with the unknown function $u(x,y,t)$: 
\begin{equation}
\label{eq:KP1}
 \frac{\partial}{\partial x}\left( 4u_t - 6uu_x -u_{xxx} \right)\, = \, 3u_{yy}.
 \end{equation}
Here $x,y$ and $t$ are the spatial and temporal coordinates respectively. The subscripts $x,y,t$ denote the partial derivatives. The function $u(x,y,t)$ displays the progression of long waves of small amplitude with slow dependence on the transverse coordinate $y$ in time $t.$ In the algebro-geometric approach to the KP equation~$\eqref{eq:KP1},$ solutions arise in the following form: 
\begin{equation}
\label{eq:u_tau}
 u(x,y,t) = 2\frac{\partial^2}{\partial x^2} \log \tau(x,y,t)+c,
\end{equation}
where $c$ is a complex constant. The function $\tau(x,y,t)$ is known as the \emph{tau function} (or \emph{$\tau$-function)} in the literature of integrable systems. A necessary and sufficient condition for a function $\tau$ to be a tau function is to satisfy the quadratic PDE, namely the Hirota bilinear equation: 
    \begin{equation}\label{eq:hirota}
( \tau_{xxxx}\tau -4\tau_{xxx}\tau_x + 3\tau_{xx}^2)
+4 (\tau_{x}\tau_t - \tau \tau_{xt}) + 6 c ( \tau_{xx}\tau -\tau_x^2)
+3 (\tau \tau_{yy} - \tau_y^2) + 8d\tau^2 
= 0.
\end{equation}

There are the \emph{quasi-periodic} solutions with tau functions~\eqref{eq:hirota} that come from complex algebraic curves~\cite{Kri1977, Dub,BelBobEnoItsMat}. In what follows, we will introduce the underlying data for such solutions in more details. Let $\mathcal{C}$ be a complex algebraic curve of genus $g,$ which is a one-dimensional complex algebraic variety as introduced in~\eqref{def_variety_complex}. Let $R$ be a \emph{Riemann matrix}~\cite{Dub} for $\mathcal{C}.$ This is a $g\times g$ complex symmetric matrix whose real part is negative definite. The \emph{Riemann theta function} associated with~$R$ is the following analytic function from $\mathbb{C}^g$ to $\mathbb{C}$:
\begin{equation}
\label{eq:RTFreal}
\theta({\bf z}\, |\, R) \coloneqq
\sum_{{\bf n} \in \mathbb{Z}^g} {\rm exp} \left( {\bf n}^t R {\bf n} + {\bf n}^t {\bf z} \right).
\end{equation}
A fundamental result of Krichever~\cite{Kri1977} asserts that there are $g$-vectors $U,V,W\in \mathbb{C}^g$ with $U\neq {\bf 0}$ such that 
\begin{equation}
\label{eq:thetatau}
\theta( \,Ux + V y + W t + D\,| \, R\,)
 \end{equation} happens to be a $\tau$-function~\eqref{eq:hirota} where $D \in \mathbb{C}^g$ is a parameter. 
In other words, the function~\eqref{eq:thetatau} satisfies~\eqref{eq:hirota}, this implies that the corresponding $u(x,y,t)$~\eqref{eq:u_tau} is a solution of the KP equation~\eqref{eq:KP1}.
One can construct this special class of tau functions from any smooth point of $\mathcal{C}.$ Set 
\begin{equation*}
    U=(u_1,\dots , u_g), \, \,\, \,\, \, V=(v_1,\dots , v_g), \, \,\, \,\, \,
    W=(w_1,\dots , w_g).
\end{equation*}
The \emph{Dubrovin threefold} $\mathcal{D}_\mathcal{C}$~\cite{CelAgoStu} of the algebraic curve $\mathcal{C}$ is formed by these triples $(U,V,W)$ in the weighted projective space $\mathbb{WP}^{3g-1},$ where $\deg(u_i)=1$ and $\deg(v_i)=2$ and $\deg(w_i)=3,$ for which there exist $c,d\in \mathbb{C}$ such that for all $D\in \mathbb{C}^g$~\eqref{eq:thetatau} satisfies~\eqref{eq:hirota}. This means that~\eqref{eq:u_tau} satisfies~\eqref{eq:KP1}. Hence $\mathcal{D}_\mathcal{C}$ is a 3-dimensional complex algebraic variety that parametrizes the solutions of the KP equation, which is guaranteed by a result of Krichever~\cite[Theorem 3.1.3]{Dub}.

\begin{figure}[h]
    \centering
    \includegraphics[width=0.35\textwidth]{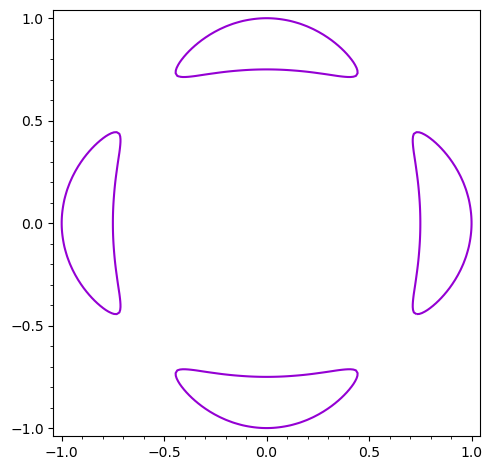}
    \includegraphics[width=0.45\textwidth]{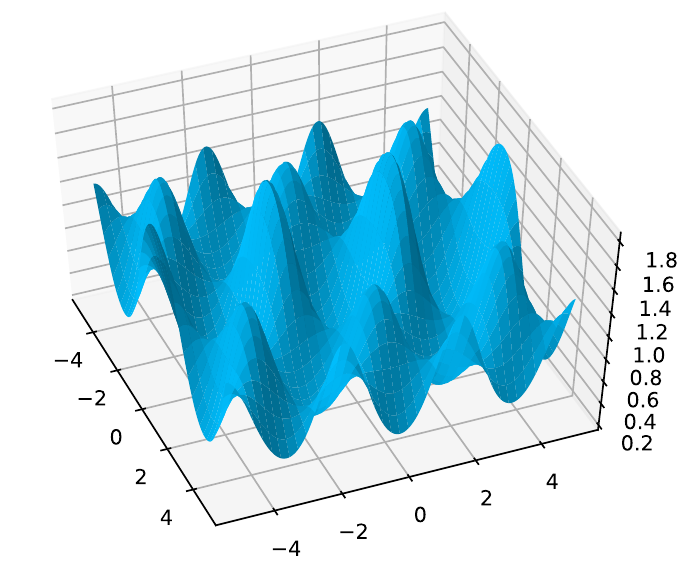}
    \vspace{-0.22cm}
    \caption{Left:~Trott curve.~Right:~The wave derived from the Trott curve whose parameters are~\eqref{eq:pointondubrovin} at $t=0$~\cite{CelAgoStu}.}
    \label{WavefromTrott}
\end{figure}

The Dubrovin threefold has a parametrization that is fully algebraic~\cite[Section 3]{CelAgoStu}. In particular, if the curve $\mathcal{C}$ is defined over the rational numbers, so is the Dubrovin threefold. This enables us to apply symbolic computational methods from nonlinear algebra using~Gr\"{o}bner bases \cite[Chapter 1]{MS2020}. In particular, we may intend to obtain the vanishing ideal of the threefold. Following~\cite[Example 3.5]{CelAgoStu}, we consider the smooth plane quartic, namely the Trott curve with its affine plane model given by the polynomial:
\begin{equation}
\label{eq:prominent}
144(x^4+y^4)-225(x^2+y^2)+350x^2y^2+81. 
 \end{equation}
The point $(0,1)$ lies on the Trott curve.~We compute a corresponding point to $(0,1)$ on the Dubrovin threefold:
\begin{equation}
\label{eq:pointondubrovin} \begin{small}
\biggl(0,  -\frac{1}{126},-\frac{1}{126},-\frac{1}{126} ,0,0,0,-\frac{1550}{55566},-\frac{1325}{37044} \biggr)
\end{small}
 \in \mathbb{WP}^8.
\end{equation}
In addition, via the {\tt SageMath}~\cite{sagemath} package \texttt{RiemannSurfaces}~\cite{BruSijZot}, we compute a Riemann matrix of the curve, we can then compute the corresponding KP solution~\eqref{eq:u_tau} by computing the corresponding Riemann theta function via the {\tt Julia} \cite{bezanson2017julia} package \texttt{Theta.jl}~\cite{AgoChu}. As the Trott curve is an {\em $M$-curve}, i.e., it has the maximum number of ovals (see \Cref{WavefromTrott}), we can proceed our computations over the real numbers~\cite{Sil}. Indeed, one can compute the Riemann matrix with real entries in this case. \Cref{WavefromTrott} illustrates the solution at $t=0$ associated with~$\eqref{eq:pointondubrovin}.$ 
It turns out that the Dubrovin threefold of the Trott curve is minimally generated by the following six polynomials---up to {saturation} by the ideal $\langle u_1,u_2,u_3 \rangle$~\cite[Chapter 3]{MS2020}: 
\begin{equation}
\label{eq:sixfortrott}
\begin{matrix}
450 u_1^2 u_3+450 u_2^2 u_3-324 u_3^3+u_2 v_1-u_1 v_2\,,\,\,
     700 u_1^2 u_2+576 u_2^3-450 u_2 u_3^2+u_3 v_1-u_1 v_3,\ \\
     576 u_1^3+700 u_1 u_2^2-450 u_1 u_3^2-u_3 v_2+u_2 v_3,  \\
     450 u_1 u_3 v_1+450 u_2 u_3 v_2+225 u_1^2 v_3+225 u_2^2 v_3-486 u_3^2 v_3+u_2 w_1-u_1 w_2, \\
     700 u_1 u_2 v_1+350 u_1^2 v_2+864 u_2^2 v_2-225 u_3^2 v_2-450 u_2 u_3 v_3+u_3 w_1-u_1 w_3, \\
     864 u_1^2 v_1+350 u_2^2 v_1-225 u_3^2 v_1+700 u_1 u_2 v_2-450 u_1 u_3 v_3-u_3 w_2+u_2 w_3,
\end{matrix}     
 \end{equation}
for which we used {\tt Macaulay2}~\cite{M2} to compute. More generally, \cite[Theorem 3.8]{CelAgoStu} gives an explicit description of the vanishing ideal of the Dubrovin threefold when $\mathcal{C}$ is a smooth plane quartic (of genus 3), for higher genus see~\cite[Theorem 5.3]{CelAgoStu}. It would be desirable to understand the geometry of the Dubrovin threefold further from the algebraic point of view. 
 
We now change the course of our discussion and review the numerical algebraic geometric side of the Dubrovin threefold: one can find polynomials over $\mathbb{C}$ that vanish on $\mathcal{D}_\mathcal{C}$ from the Hirota bilinear equation~\eqref{eq:hirota}. Namely, plugging the $\tau$-function~\eqref{eq:thetatau} into~\eqref{eq:hirota}, one obtains polynomials in $U, V, W$ and $c, d$ whose coefficients are expressions in \emph{theta constants}~\cite[Section 4]{CelAgoStu}, which were derived in~\cite[Section 4.3]{Dub}. The theta constants are values in the complex numbers which are evaluations of {theta functions with characteristic} and their partial derivatives at ${\bf 0}\in \mathbb{C}^g.$ The theta functions are translations of the Riemann theta function~\eqref{eq:RTFreal} by an exponential factor \cite[Chapter 2]{Mum1}.
Hence, we can compute polynomials defining the Dubrovin threefold $\mathcal{D}_\mathcal{C}$ via numerical evaluation of theta constants associated with $\mathcal{C}.$ Some notable mathematical software packages for computing the theta constants are \texttt{Theta.jl}~\cite{AgoChu} in {\tt Julia}, \texttt{algcurves}~\cite{algcurves} in {\tt Maple}, \texttt{abelfunctions}~\cite{abelfunctions} in {\tt SageMath}, and a package~\cite{matlabtheta} in {\tt MATLAB}, which all can be considered as tools in numerical algebraic geometry; for algebraic computations of them, see~\cite{Cel}. In fact, the theta functions and their computational aspects appear in several branches of mathematics apart from integrable systems, such as algebraic geometry~\cite{FarGruMan, AgoCelEke}, number theory~\cite{EicZag}, cryptography~\cite{Gau}, discrete mathematics~\cite{RegSte}, and statistics~\cite{AgoAme}.
\begin{figure}[h]
    \centering
    \includegraphics[width=0.43\textwidth]{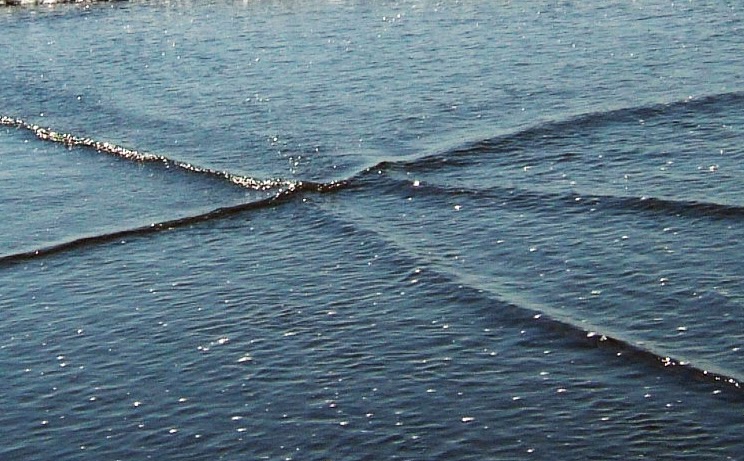} \quad
    \includegraphics[width=0.45\textwidth]{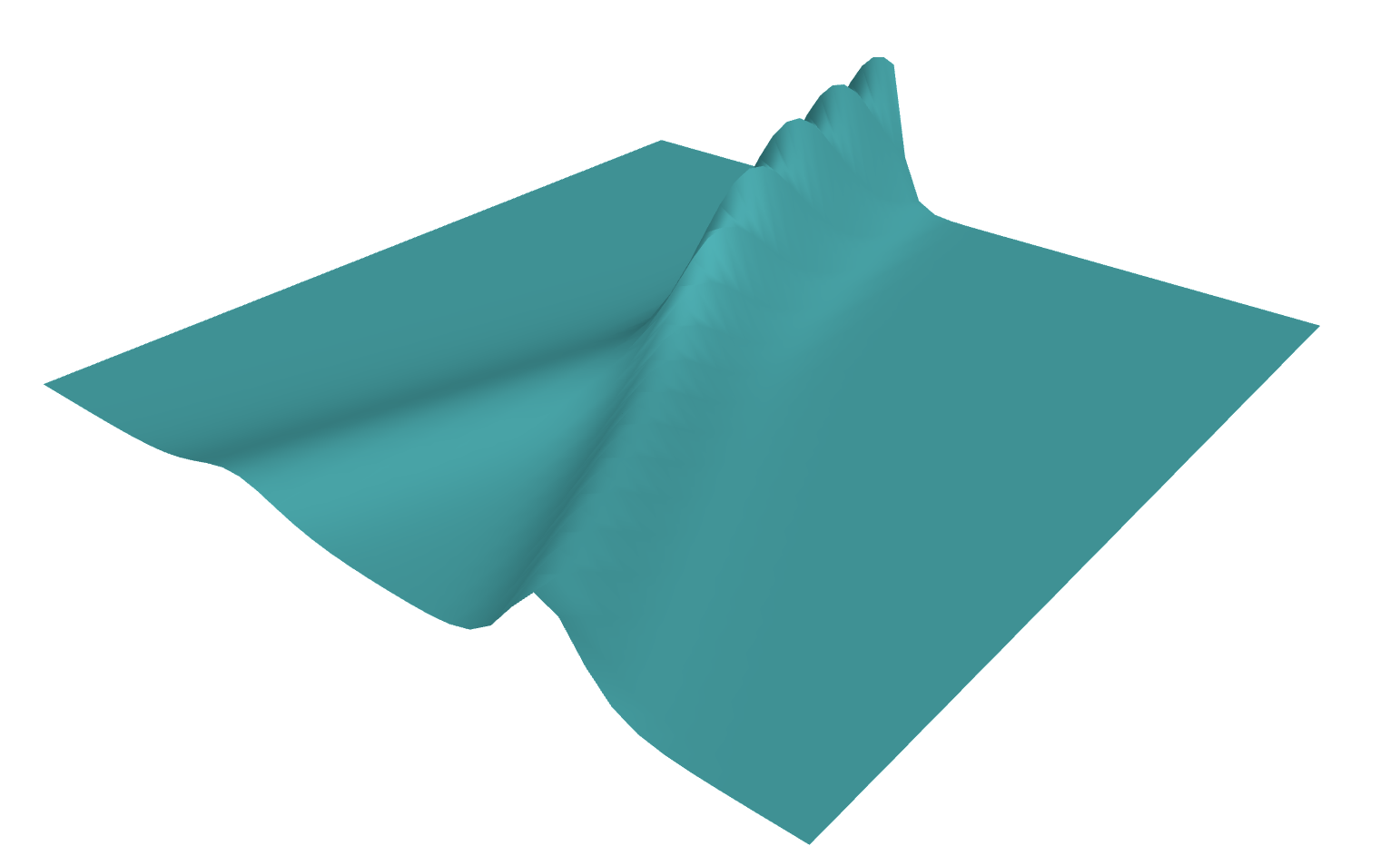}
    \caption{Left: A soliton wave that is taken in Nuevo Vallarta, Mexico by Ablowitz \cite{AblBal,Abl}. Right: A Y-soliton.}
    \label{solitonImage}
\end{figure}

 Finally, we land in combinatorial algebraic geometry with a view towards computations. A~natural question arises: what happens to the Dubrovin threefold when the underlying curve degenerates? One may expect that a potential Dubrovin variety parametrizes deformed solutions that conceivably come from singular curves. A certain class of such exact solutions comprises the so-called \emph{KP solitons}~\cite{Kod} (see~\Cref{solitonImage}), which are regular in the $xy$-plane and localized along certain rays in the plane. The soliton solutions are classified by finite-dimensional {Grassmannians}~\cite{KodWil}. There have been approaches to relate these types of solutions to singular curves~\cite{SegWil,AbeGri,Naka,KodXie}, among which the {universal Grassmannian manifold} (also known as the \emph{Sato Grassmannian}) plays a fundamental role. This infinite-dimensional Grassmannian can be viewed as a parameter space of formal power series solutions of the KP equation~\cite{Sato}. On the other hand, thinking of the quasi periodic solutions arise via the Riemann theta function, one may also study degenerations of the Riemann theta function for establishing such a relation. In fact, if one considers tropical degenerations of algebraic curves, the Riemann theta function happens to be a finite sum of exponentials \cite[Theorem 4]{AgoCelSruStu}, \cite[Theorem 3]{AgoFevManStu}. The \emph{Hirota variety} parametrizes all the tau functions arising from such a sum~\cite{AgoFevManStu}. It would be worthwhile to discover more about possible connections of the Hirota variety with the soliton solutions~\cite[Example 11]{AgoFevManStu}, also its potential link to the Dubrovin variety. In dimension 3, the Riemann theta function degenerations recover the so-called \emph{theta surfaces}~\cite{AgoCelSruStu}, classically the \emph{double translation surfaces}~\cite{Lit}. For highly singular curves, these surfaces happen to be algebraic. It turns out that the algebraic ones are given by \emph{Schur--Weierstrass polynomials}~\cite{BucEnoLey}, which gives rise to rational solutions of the KP hierarchy~\cite{AdlMos, BucEnoLey, Nak2010}. In addition, they appear in the Schur expansions of the tau functions~\eqref{eq:hirota} for solitons arising from algebraic curves~\cite[Theorem 5.1]{KodXie}. Probable relations among each of these perspectives discussed in this paragraph have not been spelled out completely yet in the literature. 
\section{Flexibility, rigidity, and configuration spaces}\label{section alex}
\emph{by Alexander Heaton}

\medskip

Configuration spaces provide interesting examples of algebraic varieties accessible and intriguing to students, practitioners, and researchers in engineering \cite{Strang-computationalscienceandengineering}, computer science \cite{SitharamStJohnSidman-HandbookGeometricConstraintSystems}, combinatorics \cite{GraverServatius-CombinatorialRigidity}, and algebraic geometry \cite{CLO}. In this section, we focus on configuration spaces arising from graphs embedded in Euclidean space. 

\begin{figure}[h!]
    \centering
    \includegraphics[width=.35\textwidth]{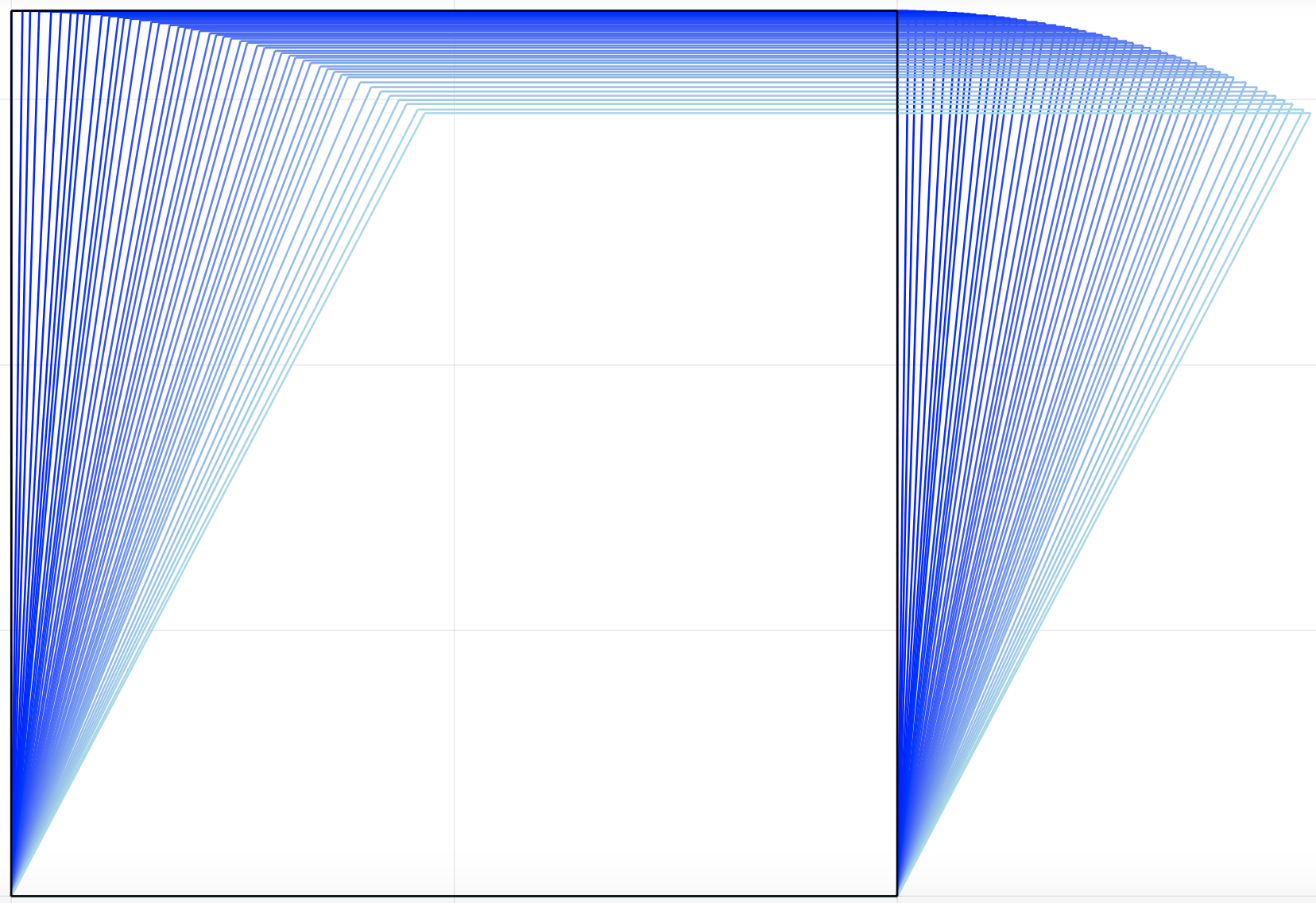}
    \includegraphics[width=.4\textwidth]{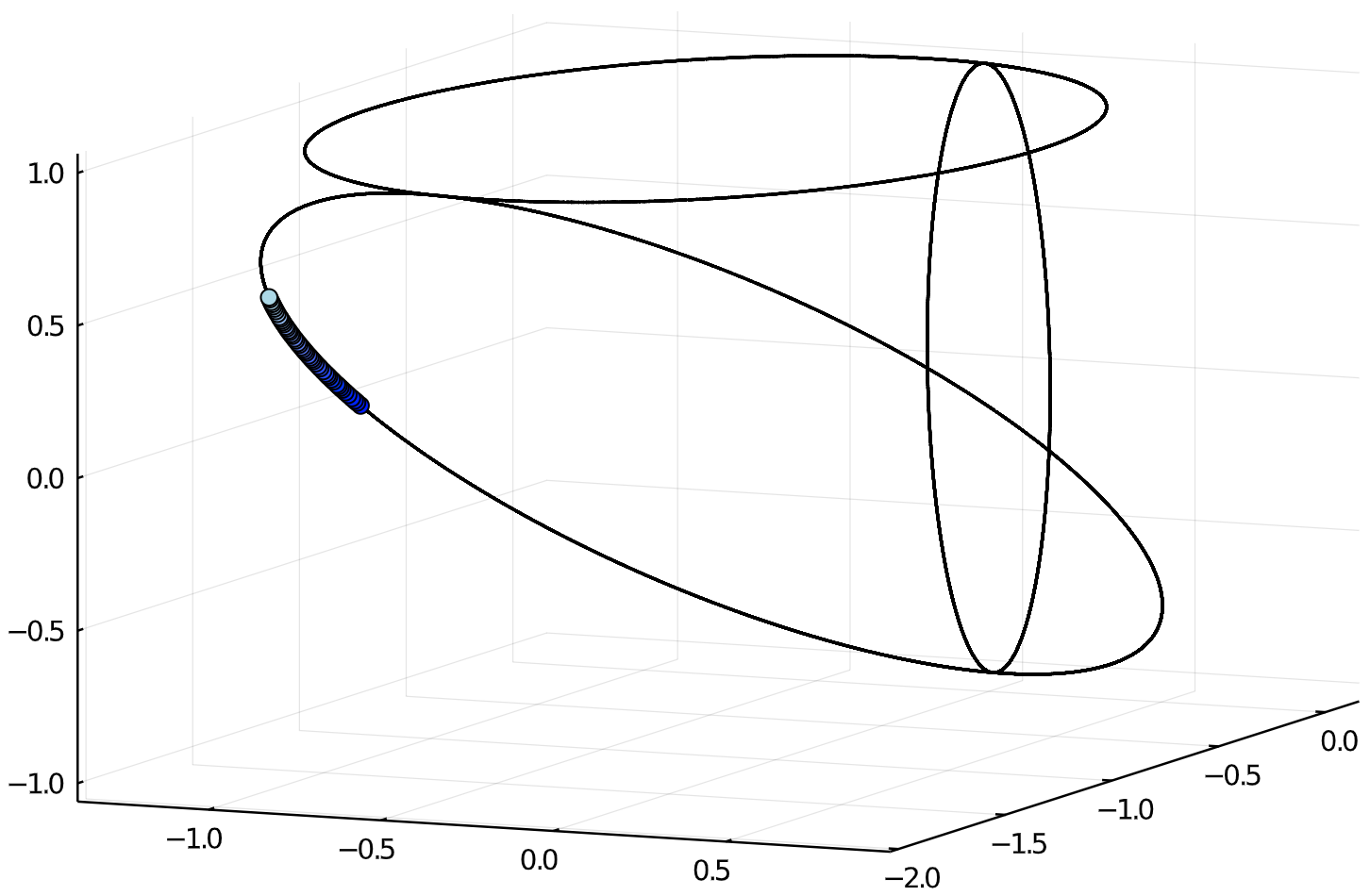}
    \caption{Configuration space of a rhombus in the plane. Left: Placements of the graph in the plane with the bottom two vertices fixed in place and the top two vertices free to move. Right: Projection of the configuration space onto a random three-dimensional subspace. The marked points on the right correspond to blue colored placements on the left, giving two ways to visualize the same data. There are only three singular points since four of the apparent intersections are artifacts of our 2d drawing of a 3d image.}
    \label{fig:rhombus-colorcoded}
\end{figure}

Before precise definitions, consider the configuration space $\mathcal{C}$ of a rhombus in the plane. A rhombus has four sides of equal length, but the angles may change. After pinning vertices $1,2$ at $(0,0),(1,0) \in \mathbb{R}^2$ and leaving the coordinates of $3,4$ as free variables, $\mathcal{C}$ is an algebraic variety of degree $6$ and dimension $1,$ with three singular points. In \Cref{fig:rhombus-colorcoded}, right, we plot its orthogonal projection onto a random three-dimensional subspace of its ambient $\mathbb{R}^4.$ Each point $p \in \mathcal{C} \subseteq\mathbb{R}^4$ corresponds bijectively to a distinct \textit{placement} $p:V \to \mathbb{R}^2,$ which is a map placing the vertices of its underlying graph $(V,E)$ with $V=\{1,2,3,4\}$ and $E=\{(1,2),(1,4),(2,3),(3,4)\}$ in the plane $\mathbb{R}^2,$ while preserving the lengths of its edges. We state the general setup below.

Let $G$ be a graph on vertices $V = \{1,2,\dots,n\}$ and edges $E \subseteq\{ \{i,j\} \mid 1\leq i < j \leq n \}.$ We will abbreviate $\{i,j\} \in E$ as $ij \in E.$ Given an edge-length function $\ell:E \to \mathbb{R}_{\geq 0},$ we are interested in {placements} (also called \textit{embeddings} or \textit{realizations}) of the graph~$G,$ which are maps $p:V \to \mathbb{R}^d$ for some dimension $d \in \mathbb{N}$ such that the Euclidean distance between $p(i)$ and $p(j)$ as points in $\mathbb{R}^d$ is equal to the value $\ell(ij).$ 
We have the polynomial map
\begin{align}\label{eq:squared-edge-length-map}
    f:\mathbb{R}^{d|V|} \to \mathbb{R}^{|E|}, \;
    \mathbf{x} \mapsto (f_{ij}(\mathbf{x}))_{ij \in E}, \quad \text{ where } f_{ij} = \sum_{k=1}^d (x_{ik} - x_{jk})^2,
\end{align}
and $\mathbf{x} = (x_{11}, x_{12}, \dots, x_{|V|d}) \in \mathbb{R}^{d|V|}.$ Therefore, the set of all placements of $G$ given $\ell$ in dimension~$d$ is the preimage of the point $(\ell(ij)^2)_{ij \in E} \in \mathbb{R}^{|E|}$ under the map $f,$ and hence has the structure of an affine algebraic variety. We call this variety the \textit{configuration space}. It is in bijection with the set of all placements satisfying the edge-length equations. For example, the rhombus in~\Cref{fig:rhombus-colorcoded} has four vertices placed in the plane $\mathbb{R}^2,$ so $\mathcal{C} \subseteq\mathbb{R}^{8}$, at first. However, by fixing the bottom two vertices we remove rigid motions and reduce the ambient dimension so that $\mathcal{C} \subseteq \mathbb{R}^4$. Next, we describe this in more detail.

Since the group of rigid motions, i.e., translations, rotations, and reflections, 
acts on the set of placements for $(G,\ell,d),$ we often remove this action by strategically {pinning vertices}. For $d=2$ this means we take $p(1) = (0,0) \in \mathbb{R}^2$ and $p(2) = (x_{21}, 0) \in \mathbb{R}^2.$ For $d=3$ we take $p(1) = (0,0,0),$ $p(2) = (x_{21}, 0, 0),$ and $p(3) = (x_{31}, x_{32}, 0).$ In fact, if we rename the vertices so that $12 \in E,$ we can also replace the variable $x_{21}$ by $\pm \ell(12).$
More generally, we can restrict the map \eqref{eq:squared-edge-length-map} %
to any subset $X \subseteq\mathbb{R}^{d|V|},$ and define the configuration space $\mathcal{C}|_X$ as the preimage of the point $(\ell(ij)^2)_{ij \in E} \in \mathbb{R}^{|E|}$ under the map $f|_X:X \to \mathbb{R}^{|E|}.$ Thus we can pin vertices to specified locations, or restrict them to move along a line or a sphere, among many possibilities.

A nice example 
in \cite{ConnellyServatiusHigher-orderRigidity} gives a framework whose configuration space has a cusp singularity. 
In this example, $n=11.$ Fixing vertices $1,6,11$ and leaving $2,3,4,5,7,8,9,10$ free, we have $\mathcal{C} \subseteq\mathbb{R}^{16}.$
\begin{figure}[h]
    \centering
    \includegraphics[width=0.4\textwidth]{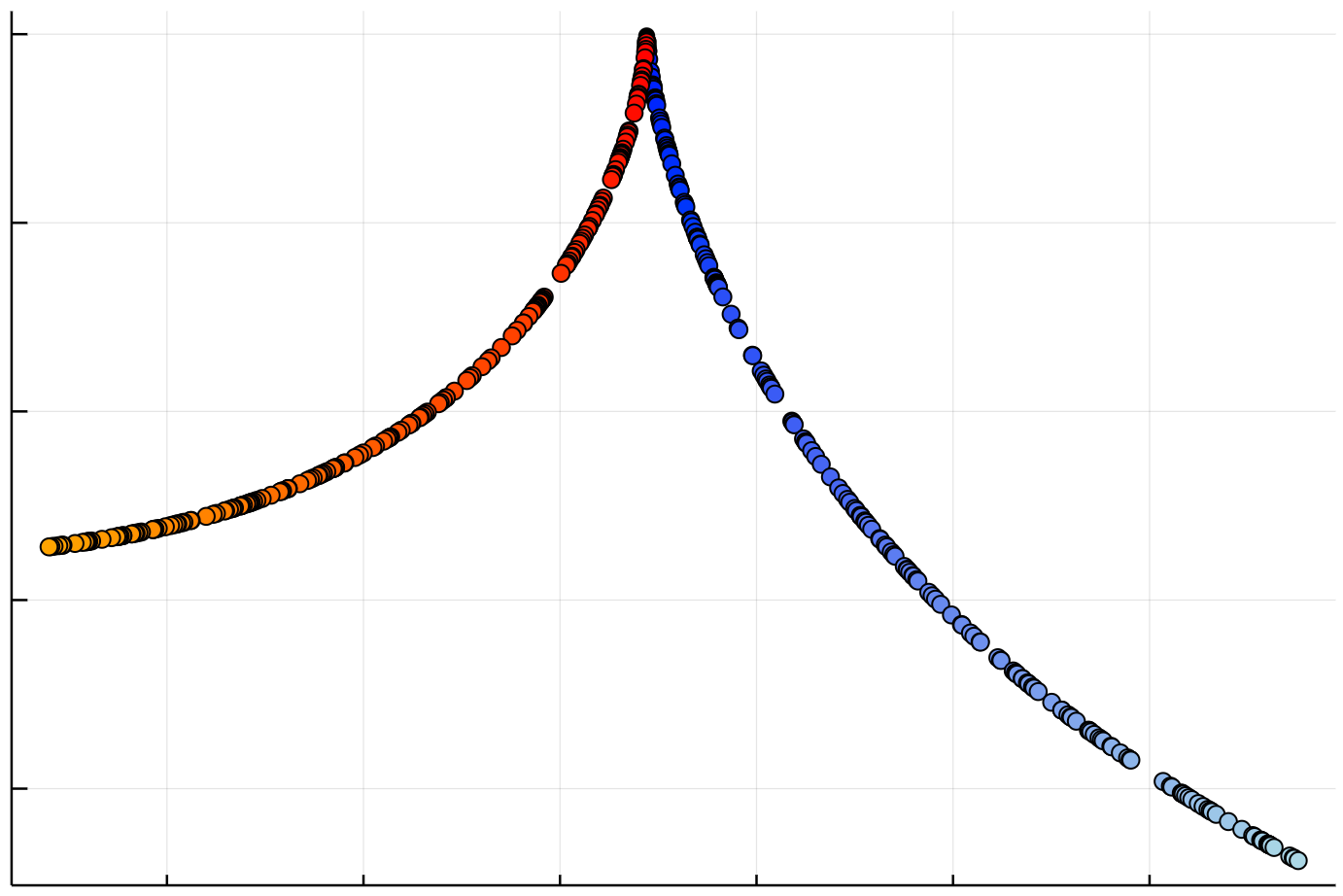}
    \includegraphics[width=0.4\textwidth]{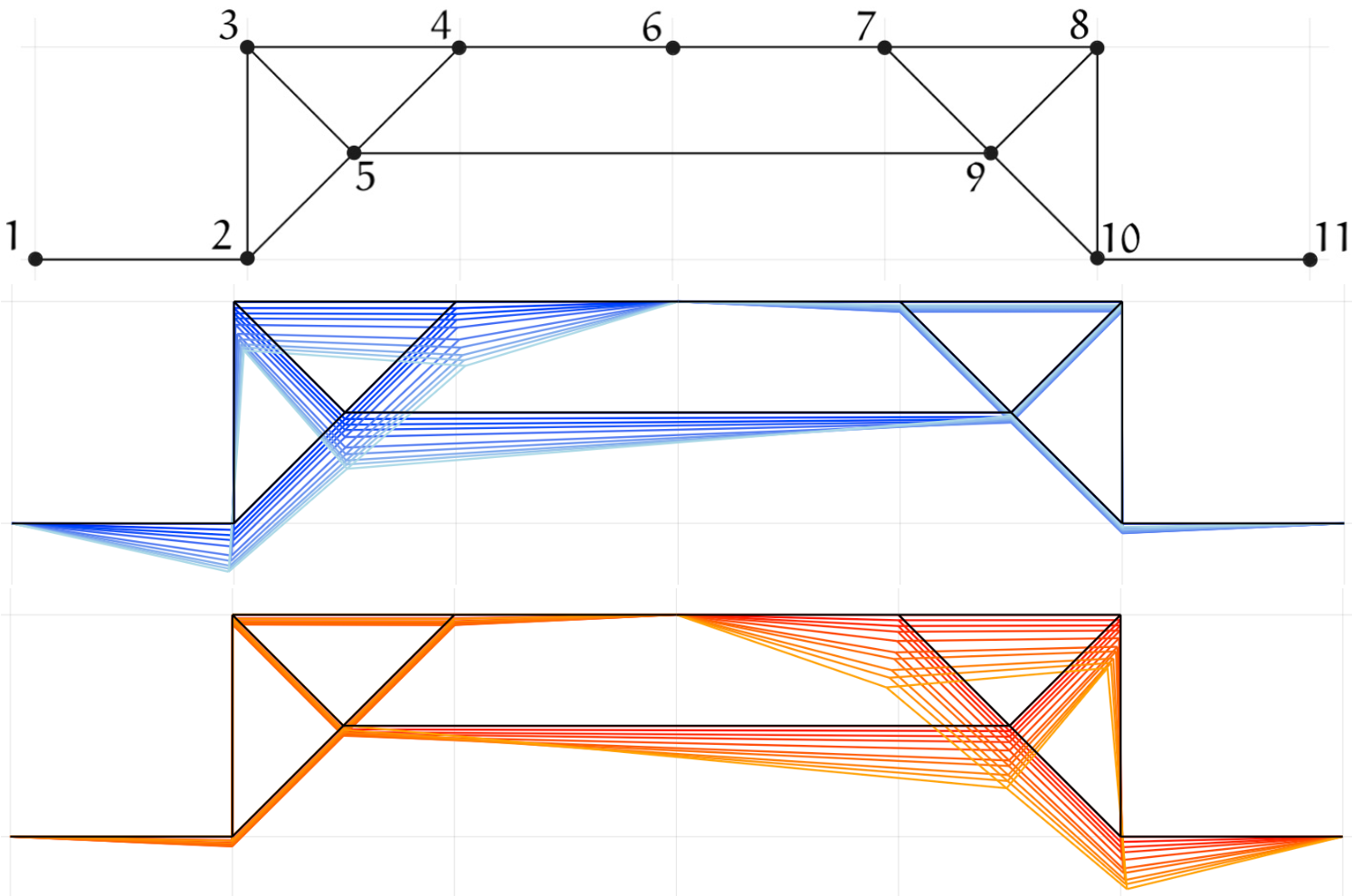}
    \caption{The left and right images are color-coded to match. Left: We project $401$ points $p^{(i)} \in \mathcal{C}$ onto a random two-dimensional subspace of $\mathbb{R}^{16}.$ $200$ orange $\to$ red points approach the singular point $p^{(i)} \to p^\star$ along one branch of the cusp, and another $200$ light-blue $\to$ blue points $p^{(j)} \to p^\star$ approach along the other branch. Right: We view each point $p^{(i)}$ as a placement map $p^{(i)}:V \to \mathbb{R}^2$ sending eleven vertices to the plane, rather than as points $p^{(i)} \in \mathbb{R}^{16}.$ Vertices $1,6,$ and $11$ are pinned and immobile. Right Top: Singular placement $p^\star.$ Right Middle: $200$ light-blue $\to$ blue placements $p^{(i)} \to p^\star$ moving toward the singular placement $p^\star$ along one branch of the cusp. Right Bottom: $200$ orange $\to$ red placements moving toward the singular placement $p^\star$ along the other branch.}
    \label{fig:cusp-with-colors}
\end{figure}
\Cref{fig:cusp-with-colors} displays a projection of one neighborhood of the configuration space from $\mathbb{R}^{16}$ to a random two-dimensional subspace, as well as the corresponding placements $p:V \to \mathbb{R}^2$ of the associated graph on eleven vertices $V.$
The colors in \Cref{fig:cusp-with-colors} match and are really two ways of visualizing the same information. By computing with homotopy continuation software \cite{HC.jl}, we discovered that $\mathcal{C}$ 
is an algebraic variety of dimension~$1$ and degree~$592.$

Among the many problems one might consider, we note the following:
\begin{enumerate}
    \item Given $(V,E)$ and an edge-length function $\ell:E \to \mathbb{R}_{\geq 0},$ understand the set of all placements  $\mathcal{C}_\ell \coloneqq \{ p:V \to \mathbb{R}^d$ : $\lVert p(i) - p(j) \rVert = \ell(ij) \text{ for all } ij \in E\}$ corresponding to $\ell.$
    \item Given $(V,E)$ and a generic placement $p:V \to \mathbb{R}^d,$ understand the set of all placements $\mathcal{C}_p \coloneqq \{ q:V\to \mathbb{R}^d : \lVert q(i) - q(j) \rVert = \lVert p(i) - p(j) \rVert \text{ for all } ij \in E \}$ corresponding to $p.$
    \item Given $(V,E)$ and a specific (non-generic) placement $p:V \to \mathbb{R}^d,$ understand the set $\mathcal{C}_p.$
\end{enumerate}

In our two examples above, the algebraic variety was one-dimensional. However, particular attention has been given to the zero-dimensional case. A placement $p:V \to \mathbb{R}^d$ is called \textit{rigid} if the local dimension of its associated configuration space is zero, and \textit{flexible} otherwise, assuming we have removed rigid motions. 
Equivalently, a placement $p$ is rigid if all nearby placements with equal edge lengths are obtained by rigid motions, and flexible if any neighborhood of $p$ contains placements which are not obtained by rigid motions, yet have the same edge lengths. Given a graph $(V,E),$ one might wonder if it has rigid placements, and if so, how many? It turns out that for \textit{generic} placements, rigidity depends only on the combinatorics of the graph. 
Originally discovered in 1927, the results in \cite{Pollaczek-Geiringer1927} and \cite{Laman1970} characterize generic 
rigidity in the plane, giving combinatorial criteria on the graph $(V,E)$ that determine whether generic placements $p:V \to \mathbb{R}^2$ will be rigid or flexible. An analogous result for graphs placed in $\mathbb{R}^3$ (or higher) remains an important and apparently difficult open problem.

Addressing a theme from the introduction, why should we use the complex numbers to study these problems, when the problem statements mention only the reals? We give two answers. First, genericity arguments using algebraic geometry are useful. For example, a placement $p:V \to \mathbb{R}^d$ is \textit{globally rigid} if every other placement with the same edge lengths is related by a rigid motion. Equivalently, global rigidity means the configuration space is zero-dimensional, consisting of a single real-valued point, after removing rigid motions. A celebrated result whose two directions are proved in \cite{GortlerHealyThurston} and \cite{Connelly-GenericGlobalRigidity} shows 
that a generic framework is globally rigid if and only if it admits a weighted graph Laplacian 
with $d+1$-dimensional kernel. The proof in \cite{GortlerHealyThurston} uses genericity arguments relying on the complex algebraic variety, including Gauss fibers and contact loci, for example.

\begin{figure}
    \centering
    \includegraphics[width=0.85\textwidth]{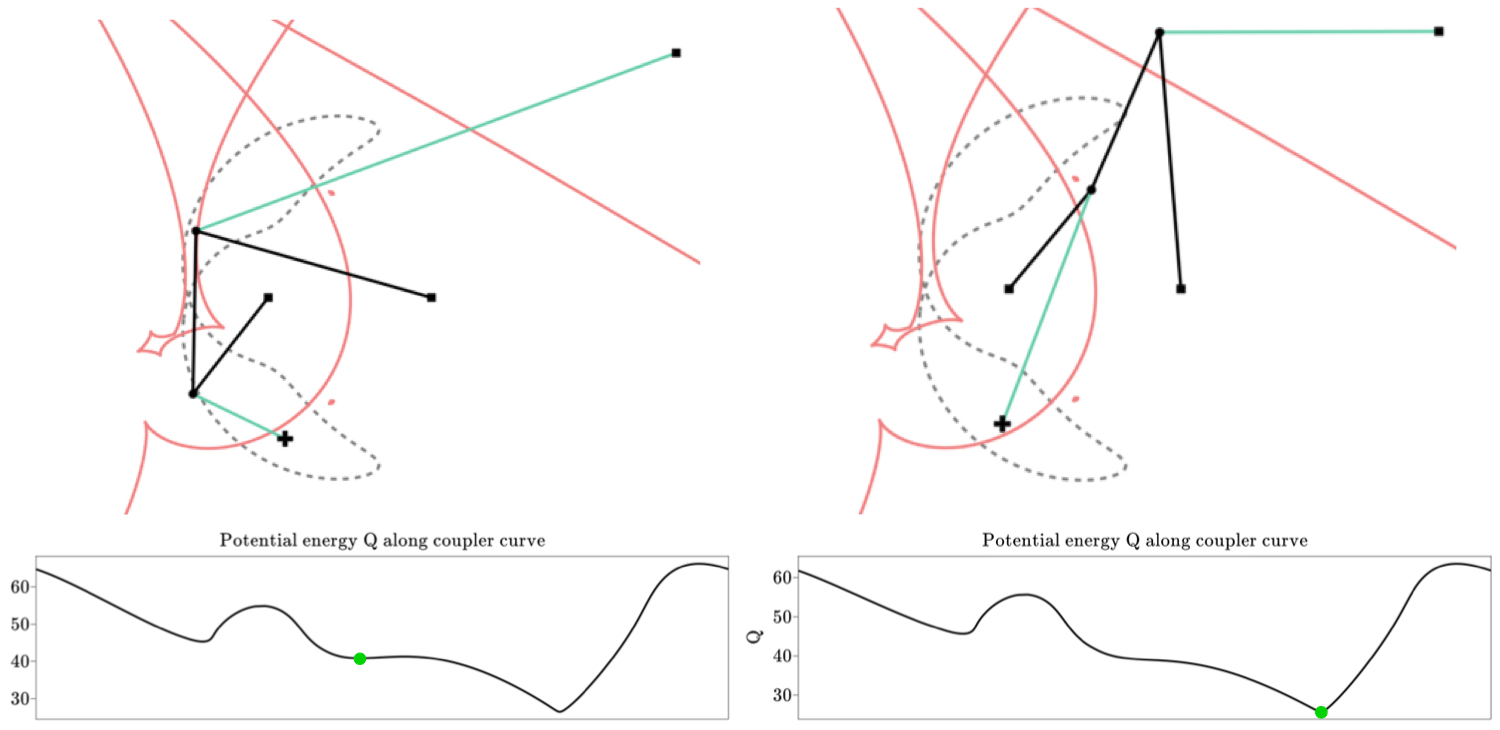}
    \caption{\footnotesize Three rigid bars in black, two elastic cables in green. Left: The elastic framework in a stable configuration. Right: Configuration of the framework after crossing the \textit{catastrophe discriminant}, depicted in red.
    The three square vertices are pinned, the cross vertex is controlled, and the two circular vertices are free: their position is found by minimizing energy over the configuration space, which is visualized by the grey, dashed \textit{coupler curve}.
    Bottom: the energy function along the coupler curve with the current position depicted in green.}
    \label{fig:fourbar_before_after}
\end{figure}

For a second reason to use complex numbers, we connect with \Cref{section polopt} on polynomial optimization. Consider a framework built of rigid bars and elastic cables. Modelling the potential energy of the cables by Hooke's law, one expects the framework to assume a position of minimum energy. The rigid bars and joints give the edges and vertices of a graph placed in space \mbox{$p:V \to \mathbb{R}^d,$} which has its associated variety $\mathcal{C}.$ The elastic cables give an energy function, and the stable positions of the framework are the placements $p \in \mathcal{C}$ which are local minima for this energy function, restricted to $\mathcal{C}.$ This is exactly the setup described in Equations \eqref{eq:problem} and \eqref{eq: KKT} of  \Cref{section polopt}. If the location of some vertices or the lengths of some bars are controlled by robotic actuators or humans, the corresponding stable local minima will also change. If the control parameters change smoothly, the stable minima will usually change smoothly as well. However, sometimes a smooth change in control parameters will cause a local minimum to disappear, perhaps merging with a local maximum or saddle. 
In these cases, a {catastrophe} occurs, since the system must move rapidly to the nearest local minimum, whose location changed discontinuously. 
In~\cite{catastrophe}, numerical nonlinear algebra is used to track all local equilibria, including those with complex-valued coordinates. This allows us to predict when local minima disappear, 
and hence predict catastrophe. Without using the complex numbers, one can easily miss important appearances and disappearances of these equilibrium positions. \Cref{fig:fourbar_before_after} demonstrates a catastrophic jump in positions of a  
framework 
augmented with two green elastic cables.
\section{Biochemical reaction networks} 
\emph{by O\u{g}uzhan Y\"ur\"uk}
\medskip

Chemical reaction networks theory (CRNT) studies the behavior of chemical systems by modeling the concentrations of the species, i.e., the chemical substances that take part in the reaction network, over time.
In this section, we present a refined introduction to CRNT and address a specific problem using tools from nonlinear algebra. We refer the reader to \cite[Part 1, Chapter~3]{feinberg_2019} for a comprehensive introduction to CRNT. 
\begin{definition}
	A \emph{chemical reaction network} $\mathcal{N} = (\chi, \mathcal{R})$ consists of a set \mbox{$ \chi \coloneqq \left\{X_1, \ldots ,X_n\right\}$}, called \emph{species}, and a set of \emph{reactions} ${\mathcal{R}} \coloneqq   \left\{R_1, \ldots ,R_l\right\}$ of the following form:
	\begin{align}
	\label{Equation:GeneralReaction}
	R_j :  \sum_{i = 1}^n a_{ij} X_i \ce{->[{\kappa}_j]} \sum_{i = 1}^n b_{ij} X_i,
	\end{align}
	where $a_{ij},b_{ij} \in \mathbb{Z}_+$ are called \emph{stoichiometric coefficients} and each ${\kappa}_j \in \mathbb{R}_{\geq 0}$ is a positive real parameter called \emph{the reaction rate constant} of $R_j$ for $j= 1, \ldots , l.$
\end{definition}
A chemical reaction is said to be at a \emph{steady state} if the amounts of the species in the reaction remain constant. \emph{Multistationarity}, that is, the existence of multiple steady states, is an important notion for biochemistry due to its link to the cellular decision making and on/off-responses to graded input \cite{laurent1999,ozbudak2004,Xiong:2003jt}. 

In the context of CRNT,
there are numerous methods to detect multistationarity for some fixed choice of parameter values \cite{Feinbergss,feliu_newinj,wiuf-feliu,PerezMillan,craciun2008,conradi-PNAS,control,crnttoolbox,dickenstein2019multistationarity,dickenstein:regions}; also see \cite{joshi2015survey} for a comprehensive survey. However, determining the exact region of multistationarity in the parameter space is a difficult problem with complicated answers. Recently, in \cite{FeliuPlos}, the authors point out a class of networks whose multistationarity can be decided by studying the nonnegativity of a relevant polynomial with parametrized coefficients. Polynomial nonnegativity is not only studied classically (cf. \cite{Hilbert:1888,Hurwitz:AMGM}), but also is a significant notion that underlies polynomial optimization as in \Cref{section polopt}. Therefore, it makes sense to apply expertise around the polynomial nonnegativity to study multistationarity in CRNT. In this section, we first explain the relationship between multistationarity and nonnegativity and then give a short overview of the two methods that were used in \cite{MyCRN1} to study the multistationarity in dual phosphorylation/dephosphorylation cycle.   

Given a reaction network $\mathcal{N} = (\chi, \mathcal{R}),$ we denote the vector of concentrations \linebreak with \mbox{$\mathbf{x} = \left(x_1, \ldots, x_n\right) \in \mathbb{R}_{\geq 0},$} where each $x_i$ corresponds to the concentration of the species $X_i.$ This vector is time-dependent, hence we sometimes write $\mathbf{x}(t)$ with $t\geq 0$ to stress this dependence, and denote its first time-derivative with $\dot{\mathbf{x}}.$ Under the assumption of mass action kinetics, the network~$\mathcal{N}$ satisfies the following system of ordinary differential equations:
\begin{align}
\dot{\mathbf{x}} =  N \cdot \mathbf{v(\mathbf{x})},
\label{Equation:CRN-ODE}
\end{align}
where $N$ is an $n \times l$-matrix with entries $N_{ij} \coloneqq b_{ij}-a_{ij}$ and $\mathbf{v({x})}$ is the column vector of length~$l$ whose entries are given as $v_j(\mathbf{x}) = {\kappa}_j x_1^{a_{1j}} \cdots x_n^{a_{nj}}$ for some choice of parameter ${\kappa}_j \in \mathbb{R}_+.$
If \mbox{$\rank(N)=r$}, then the \emph{stoichiometric subspace $S$} of the network~$\mathcal{N}$ is an $r$-dimensional linear subspace of $\mathbb{R}^n$ generated by the columns of $N.$
For a given initial vector $\mathbf{x}(0) \in \mathbb{R}_{\geq 0}^n$ of concentrations,  $\mathbf{x}(t)$ is confined in the affine space $\left( \mathbf{x}(0) + S \right) \cap \mathbb{R}_{\geq 0}^n$ for $t\geq 0$ by \cite[Lemma 3.4.5]{feinberg_2019}.
Let $W$ be a full-rank row reduced $(n-r) \times n$-matrix such that $WN=0,$ and define $\mathbf{c}\coloneqq W \cdot \mathbf{x}(0) \in \mathbb{R}^{n-r}.$ Then, one can identify the affine space $ \mathbf{x}(0) + S$ with the vector $\mathbf{c} \in \mathbb{R}^{n-r},$ because $W\mathbf{x} = \mathbf{c}$ for any $\mathbf{x} \in \left( \mathbf{x}(0) + S \right).$

If a point $\mathbf{x^*} \in \mathbb{R}^n$ satisfies that $ N \cdot \mathbf{v(\mathbf{x^*})} = {0},$ then it is called a \emph{steady state} of the network~$\mathcal{N}.$ If we fix the parameter vector $ \mathbf{\kappa} \coloneqq \left( \kappa_1, \ldots, \kappa_l \right) \in \mathbb{R}_+^l,$ then $V\coloneqq \left\{\mathbf{x} \in \mathbb{R}^n \mid \dot{\mathbf{x}} =  N \cdot {v(\mathbf{x})}= {0} \right\}$ is an algebraic variety in $\mathbb{R}^n$ given as the zero set of $n$ polynomials in $\mathbb{R}[x_1,\ldots,x_n].$
In particular, it is an intriguing question to find the number of equilibrium points confined in the the affine space $ \mathcal{P}_{\mathbf{c}}\coloneqq \left\{ \mathbf{x} \in \mathbb{R}_{\geq 0}^n \mid W\mathbf{x} = \mathbf{c} \right\},$ i.e., the cardinality of
\begin{equation*}
    V \cap \mathcal{P}_{\mathbf{c}} = \left\{\mathbf{x}\in \mathbb{R}_{\geq 0}^n \mid N \mathbf{v(x)} = 0 \right\} \cap  \left\{\mathbf{x}\in \mathbb{R}_{\geq 0}^n \mid W\mathbf{x} = \mathbf{c} \right\}.
\end{equation*}
Given a stoichiometric compatibility class $\mathcal{P}_{\mathbf{c}} $, it is an intriguing question to find the positive steady states within this class. Therefore, we define $\mathcal{P}_{\mathbf{c}}^+ :=  \left\{\mathbf{x}\in \mathbb{R}_{> 0}^n \mid W\mathbf{x} = \mathbf{c} \right\}$, and name $\mathcal{P}_{\mathbf{c}}^+$ as  \emph{the postive stoichiometric compatibility class associated to $\mathcal{P}_{\mathbf{c}}$}.

If there exists a compatibility class $\mathcal{P}_{\mathbf{c}}$ that contains more than one steady state for some~$\mathbf{\kappa} \in \mathbb{R}_+^l,$ then the parameter $\mathbf{\kappa}$ is said to \emph{enable multistationarity} in the network $\mathcal{N}$. Deciding the existence of parameters that enable multistationarity for a given reaction network is a central question in CRNT. Using the particular approach that was introduced in \cite{FeliuPlos}, this question reduces down to checking the sign of the determinant of a relevant Jacobian matrix. This method in particular works in the case when $\mathcal{P}_{\mathbf{c}}^+$ admits a positive parameterization and the network is \emph{conservative}, i.e., if each $\mathcal{P}_{\mathbf{c}}$ is a compact set (see \cite[Supplementary Info Section 3.2]{FeliuPlos}). A positive parameterization of the set $V \cap \mathcal{P}_{\mathbf{c}}^+$ is a surjective function $\phi:\mathbb{R}^m \to \mathbb{V \cap \mathcal{P}_{\mathbf{c}}^+}$ for some $m<n$. To construct this Jacobian matrix, first let us denote the index of the first nonzero entry for each row of $W$ with $i_1,\ldots,i_{n-r}.$ Then, we consider the map $\varphi_c: \mathbb{R}_{\geq 0}^n \to \mathbb{R}^n$ whose $i$-th entry $\varphi_c(x)_i$ is equal to $(W\mathbf{x}-\mathbf{c})_i$ if $i = i_1,\ldots,i_{n-r},$ and it is equal to $N \mathbf{v(x)}_i$ otherwise. We note that the steady states in the stoichiometric compatibility class $\mathcal{P}_{\mathbf{c}}$ is $V \cap \mathcal{P}_{\mathbf{c}} = \left\{\mathbf{x}\in \mathbb{R}_{\geq 0}^n \mid \varphi_c(\mathbf{x}) = 0 \right\}.$

\begin{theorem}[{\cite[Theorem 1]{FeliuPlos}}]
	\label{Theorem:MainMultistationaryTheorem}
	Let $\mathcal{N} = \left(\chi,  \allowbreak  \mathbb{R}  \right)$ be a conservative chemical reaction network with stoichiometric matrix $N \in \mathbb{R}^{n \times l}$ of rank $r.$ Furthermore, let $\mathcal{P}_{\mathbf{c}}$ be a nonempty stoichiometric compatibility class without any boundary equilibrium point where $\mathbf{c} \in \mathbb{R}^{n-r}$ and let $\phi:\mathbb{R}^m \to \mathbb{V \cap \mathcal{P}_{\mathbf{c}}^+}$ be a positive parameterization of $V \cap \mathcal{P}_{\mathbf{c}}^+$ for some $m<n$. If we denote the Jacobian matrix of $\varphi_c(\mathbf{x})$ after the positive reparameterization with $M(\mathbf{x})$, then the following statements hold:
	\begin{itemize}
		\item[(a)] If $ \sign(\det(M(\mathbf{x}))) = (-1)^r $ for all $\mathbf{x} \in V \cap \mathcal{P}_{\mathbf{c}}^+,$ then there is exactly one positive equilibrium in $\mathcal{P}_{{c}}.$ 
		\item[(b)] If $ \sign(\det(M(\mathbf{x}))) = (-1)^{r+1} $ for some $\mathbf{x} \in V \cap \mathcal{P}_{\mathbf{c}}^+,$ then there are at least two positive equilibria in $\mathcal{P}_{\mathbf{c}}.$
	\end{itemize}
\end{theorem}
We will now see an example of this theorem in action. Consider the \emph{hybrid histidine kinase system}, which models the dual phosphorylation of histidine kinase enzyme. 
\begin{align*}
HK_{00} \ce{->[\kappa_1]} HK_{p0} \ce{->[\kappa_2]} HK_{0p} \ce{->[\kappa_3]} HK_{pp}, & \qquad HK_{0p} + RR \ce{->[\kappa_4]} HK_{00} + RR_p,\\
HK_{pp} + RR \ce{->[\kappa_5]} HK_{p0} + RR_p, &\qquad 
RR_p \ce{->[\kappa_6]} HK_{pp}.
\end{align*}
\noindent
In the network given above, $HK_{p0}$ (and $HK_{0p}$) denotes the histidine kinase with its first (and the second) site phosphorylated, $HK_{00}$ and $HK_{pp}$ denote the histidine kinase with zero and two phosphorylated sites, and  $RR$ denotes the response regulator protein which may be phosporylated by histidine kinase to form $RR_p.$ We note that the stoichiometric matrix $N$ is of rank 4 in this case.
If we let $X_1 = HK_{00}, \, X_2 = HK_{p0},\, X_3 = HK_{0p}, \,X_4 = HK_{pp},$ $X_5 = RR, \, X_6= RR_p,$ then 
\begin{align*}
    \varphi_c({x}) =  \big( & x_1+x_2+x_3+x_4 - c_1, \ \kappa_1x_1 - \kappa_2x_2 + \kappa_5x_4x_5, \   \kappa_2x_2 - \kappa_3x_3 - \kappa_4x_3x_5,  \\
       & x_5+x_6 - c_2, \  \kappa_4x_3x_5+\kappa_5x_4x_5-\kappa_6x_6 \big),
\end{align*}
and the determinant of the Jacobian matrix of $\varphi_c(\mathbf{x})$ is
\begin{align*}
M({x}) =\ & \kappa_2\kappa_4\kappa_5(\kappa_1-\kappa_3)x_3x_5 + \kappa_1\kappa_2\kappa_4\kappa_5x_4x_5 + \kappa_4\kappa_5\kappa_6(\kappa_1+\kappa_2)x_5^2 \\ &+ \kappa_1\kappa_2\kappa_3\kappa_4x_3 + \kappa_1\kappa_2\kappa_3\kappa_5x_4 + \kappa_1\kappa_5\kappa_6(\kappa_3+\kappa_2)x_5 + \kappa_1\kappa_2\kappa_3\kappa_6 .
\end{align*}
If $\kappa_1 \geq \kappa_3,$ then \cref{Theorem:MainMultistationaryTheorem} implies that the system has a unique equilibrium point in each stoichiometric compatibility class. For $\kappa_1 < \kappa_3,$ if we consider the positive parameterization  
\begin{align*}
    \Phi({x_4,x_5}) =  \big( \frac{\kappa_4\kappa_4x_4x_5^2}{\kappa_1\kappa_3}, \ \frac{\kappa_5(\kappa_4x_5+\kappa_3)x_4x_5}{\kappa_2\kappa_3} , \  \frac{\kappa_5x_4x_5}{\kappa_3},  \ x_4, \ x_5, \  \frac{\kappa_5(\kappa_4x_5+\kappa_3)x_4x_5}{\kappa_3\kappa_6} \big)
\end{align*}
of $V \cap \mathbb{R}^6_{>0}$, then the coefficient of the term $x_4x_5$ in the polynomial $M(\Phi({x_4,x_5}))$ is negative. Since this term corresponds to a vertex of the Newton polytope of $M(\Phi({x_4,x_5})),$ the polynomial $M(\mathbf{x})$ can take negative values (see e.g., \cite[Proposition 2.3]{MyCRN1}). Therefore, \cref{Theorem:MainMultistationaryTheorem} implies that the system enables multistationarity if $\kappa_1 <\kappa_3.$

Note that the determinant $M(\mathbf{x})$ is a polynomial in $\mathbb{R}[x_1,\ldots,x_n]$ whose coefficients are \linebreak parametrized by the entries of $\kappa.$ Unlike in the previous example, it is not always easy to identify the pa\-ra\-me\-ter region that guarantees the nonnegativity of $M(\mathbf{x}).$ There are various tools in real algebraic geometry to address this problem, and here we mention two of these approaches. In particular, we point out how these methods have been utilized for the case of 2-site phosphorylation cycle in~\cite{MyCRN1}. 

The first approach is based on quantifier elimination and the Cylindrical Algebraic Decomposition (CAD) algorithm introduced by Collins \cite{collins1975quantifier}. For an overview of quantifier elimination and CAD, we refer to \cite{caviness2012quantifier}. Via CAD, one can decompose $\mathbb{R}^n$ into connected semialgebraic sets (see \cref{section polopt} for a definition) on which $M(\mathbf{x})$ has constant sign, and determine the sign of $M(\mathbf{x})$ at each of these components. The CAD algorithm is not viable for large examples, but one can still make use of CAD by reducing the number of variables or parameters via some algebraic tricks such as restricting the polynomial of question to a smaller one with \cite[Proposition 2.3]{MyCRN1}. For instance in \cite[Section 3.1]{MyCRN1}, the authors point out a sufficient condition for nonnegativity following such a strategy and using the package {\tt RegularChains} in {\tt Maple} \cite{xiao:regularchain}.

The second approach considers symbolic nonnegativity certificates for real polynomials. There are several methods to check nonnegativity of a polynomial with explicit coefficients \cite{Blekherman:Parrilo:Thomas:SDOptAndConvAlgGeo,Laserre:MomentsPosPoly,Marshall:Book}. In fact, this is a feasibility problem for polynomial optimization. Here, we mention a particular method based on circuit polynomials, which are a special class of polynomials in $\mathbb{R}[x_1,\ldots,x_n]$ whose set of exponents forms a minimal affine dependent set in $\mathbb{R}^n.$ The circuit polynomials have been introduced in \cite{Iliman:deWolff:Circuits} and their nonnegativity is fully characterized by a symbolic condition given in terms of their coefficients and the combinatorics of its support, see \cite[Theorem 1.1]{Iliman:deWolff:Circuits}. If we decompose $M(\mathbf{x})$ into circuit polynomials, then assuming the nonnegativity of each circuit polynomial in this decomposition is sufficient to make $M(\mathbf{x})$ nonnegative as well. The nonnegativity of each circuit polynomial in the decomposition yields a symbolic condition in terms of the parameters $\kappa_1, \ldots, \kappa_l,$ and hence describes a region in the parameter space. For instance, this is used in \cite[Theorem 3.5]{MyCRN1}, where the authors formulate a sufficient condition for nonnegativity by writing $M(\mathbf{x})$ as sum of four circuit polynomials.

To conclude, using techniques from real algebraic geometry to study the signs of a parametric multivariate polynomial on the positive orthant as a function of the parameters, it is possible to find open sets in the space of parameters $\mathbb{R}_+^l$ that enable multistationarity. The techniques we mention here do not solve the problem of multistationarity. However, they vastly extend our understanding of the region of multistationarity, and are not exclusive to the dual phosphorylation cycle. Moreover, one can further apply these techniques to other problems in chemical reaction networks theory with similar flavor. For example, the study of signs plays a key role when analyzing the stability of steady states or the presence of Hopf bifurcations via the Routh--Hurwitz criterion (see, e.g., \cite{torres:stability,shiu:hopf}). 

\section{Algebraic vision}\label{section algvision}
\emph{by Timothy Duff}
\medskip 

Projective space, algebraic varieties, rational maps, and many other notions from algebraic geometry appear naturally in the study of image formation with respect to various different camera models. 
These basic notions play a distinguished role in computer vision, where one of the primary goals is to build systems capable of reconstructing 3D geometry (scene and cameras) from data in several images.
For a broad overview of computer vision, we refer to the text by Szeliski~\cite{Szeliski}.

\begin{figure}[h]
\begin{center}
\def\t{0.28}
\def\s{1.8}
\begin{tikzpicture}[scale = 0.8]
\coordinate (c) at (0,0,0);
\coordinate (p) at (-1,3/2,0);
\coordinate (x0) at (\t,-\t*3/2,0);
\coordinate (x) at (1,-3/2,0);
\coordinate (x1) at (\s,-\s*3/2,0);
\coordinate (e1) at (1,0,0);
\draw[thick, ->] (c) -- (e1) ;

\draw[fill=black!10] (1,-2,-2) -- (1,-2,2) -- (1,2,2) -- (1,2,-2) node[right] {$H$} -- cycle;
\draw[thick, dashed, -] (p) -- (x0) ;
\draw[thick, dashed, ->] (x) -- (x1) ;
\draw[fill] (p) circle (0.05) node[below left] {$(x,y,z)$};
\draw[fill] (x) circle (0.05) node[below right] {$\;\;(x/z,y/z,1)$};
\draw[fill] (e1) circle (0.05) node[above] {$(0,0,1)$};
\draw[fill] (0,0,0) circle (0.05) node[left] {$\mathbf{0}$};

\draw[thick, ->] (c) -- (0,0,2) ;
\draw[thick, ->] (c) -- (0,2,0) ;
\draw[thick, ->] (e1) -- (2,0,0) ;
\end{tikzpicture}
\end{center}
\caption{A pinhole camera with principal point equal to $(0, \, 0, \, 1) \in H$ and focal length $1.$ The point $(x,y,z)$ is projected onto the plane $H.$ The resulting image is the point $(x/z,y/z,1).$ The dashed line corresponds to the point in $\PP_{\mathbb R}^2$ which is represented in homogeneous coordinates by $[x:y:z].$}
  \label{fig:pinhole}
  \end{figure}
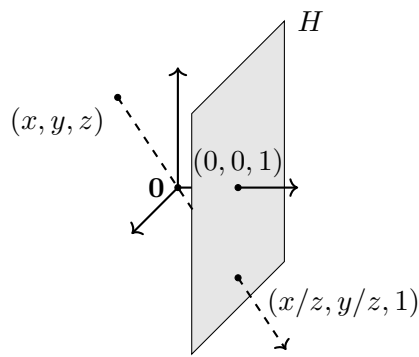

To understand the importance of geometry in computer vision, the {\em pinhole camera model} is an excellent starting point. This model is depicted in Figure~\ref{fig:pinhole}.
Here, a real-life camera is modeled as a projection from a point in space onto a plane.
The center of projection is an idealized lens, through which rays of light pass to form an image.
The coordinates of our {camera frame} are chosen so that the center of projection is the origin and so that the {image plane} is given\linebreak  by~$H=\{ (x,y,z) \in \RR^3 \mid z=1\}.$
We remark that our choice of coordinates implies that two world points $(x,y,z), (-x,-y,-z) \in \RR^3$ which differ only by sign  will produce the same image in $H.$
As a matter of convenience, it is more standard to draw the image plane in {\em front} of the camera center, unlike in our picture.

The equations of a world-to-image map for this camera, which can be derived algebraically or using similar triangles, are given by
\begin{equation}
\label{eq:R3R2}
    \begin{split}
\RR^3 \dashrightarrow H, \quad
        (x,y,z) \mapsto (x/z, y/z ,1);
    \end{split}
\end{equation}
see, e.g., \cite[Section 6.1]{HZ-2003}.
This map is nonlinear, and undefined when $z=0.$
However, it can be better understood through the lens of projective geometry.
In this approach, thoroughly laid out in the text by Hartley and Zisserman~\cite{HZ-2003}, each point in the image is naturally identified with the light-ray that passes through it---a line through the origin in $\RR^3,$ like the dashed line in Figure~\ref{fig:pinhole}.
There are also exceptional light-rays where the map~\eqref{eq:R3R2} is undefined.
These exceptional lines in space correspond to ``vanishing points" where two parallel lines in our image meet.
The \emph{real projective space} $\mathbb{P}_{\mathbb R}^2$ is the space of {all} lines through the origin in $\RR^3$ (similar as the definition of projective space in \Cref{section intro}).
Using {homogeneous coordinates} on $\PP_{\mathbb R}^2,$ we may rewrite our image coordinates as $[x/z:y/z:1] = [x:y:z].$
In doing so, we obtain a projective-linear map
\begin{equation}
\label{eq:P3P2}
    \begin{split}
        \PP_{\mathbb R}^3 \dashrightarrow \PP_{\mathbb R}^2,\quad
        [x:y:z:w] \mapsto [x:y:z],
    \end{split}
\end{equation}
which is defined for all world points except the camera center $[0:0:0:1].$

In various applications, the information provided by \emph{multiple} cameras must be combined.
These applications may involve very intensive computation (e.g.,~building 3D models of cities from large collections of photos~\cite{Rome}) or have strict real-time requirements (e.g.,~multi-camera systems for~$360^\circ$ field-of-view on an autonomous vehicle~\cite{pollefeys-multicamera}). 
A common thread throughout these applications is that the data in images are noisy and may even suffer from a large presence of gross errors, also known as \emph{outliers}.
For example, a common task is to estimate the relative orientation between two cameras when given many pairs of corresponding points $(\xx_1, \yy_1), \ldots , (\xx_m, \yy_m) \in {\PP_{\mathbb R}^2} \times {\PP_{\mathbb R}^2}$ between two images.
For two cameras whose internal parameters like focal length are known, this relative orientation is encoded by the \emph{essential matrix} ${E}$; see \cite[Section 9.6]{HZ-2003}. This is a $3\times 3$-matrix with real entries that satisfies the following polynomial constraints. 
Thinking of each ${\xx_i}$ and ${\yy_i}$ as a $3\times 1$-matrix, we must have
\begin{equation}
    \label{eq:5pp}
    \begin{split}
    2 {E} {E}^\top {E} - \tr ({E} {E}^\top ) {E} = 0_{3\times 3},\quad
    \det {E} = 0, \quad 
    {\yy_i^\top} {E} {\xx_i} =0 \,\,\,\textrm{for}\,\,\, i=1,\ldots , m.
    \end{split}
\end{equation}
For this problem, outliers are mismatched pairs $({\xx_i}, {\yy_i})$ which do not correspond to a common point seen by the two cameras.
The prevalence of outliers calls for robust estimation techniques, among which heuristics based on \emph{RANSAC} (RANdom SAmpling and Consensus~\cite{RANSAC}) have been ubiquitous in computer vision.
The key to RANSAC is using the minimal number of measurements---for the equations in~\eqref{eq:5pp}, this means $m=5$---to generate several ``hypotheses" for ${E},$ which are then checked against the remaining data.
An ingenious scheme for solving the $m=5$ case of~\eqref{eq:5pp} was proposed and implemented by Nist\'{e}r~\cite{Nister}. 
As in previous sections, this solution exploits genericity.

\begin{proposition}[{\cite[Th\'{e}or\`{e}me 6.4]{Demazure}}]
\label{prop:5pt}
For $m=5$ and generic data $(\xx_1 , \yy_1) , \, \ldots  , \, (\xx_5, \yy_5),$ there are exactly $10$ complex matrices $E$ up to scale (i.e.~in $\PP^8$) which satisfy \eqref{eq:5pp}.
\end{proposition}

Nist\'{e}r's five-point algorithm has since seen been incorporated into state-of-the-art 3D reconstruction pipelines, in which typically the initial estimates obtained through RANSAC are further refined with nonlinear least-squares. 
A sizeable literature has also emerged on solving other minimal problems in real-time with specialized solvers based on Gr\"{o}bner bases and resultants (see, e.g.,~\cite{AutoMin,syz}).
Recently, homotopy continuation methods have played a complementary role, enabling the study of a much wider class of minimal problems than \mbox{before~\cite{Kil17,TRPLP,PLMP}.}

So far, we have barely scratched the surface of how nonlinear algebra can be applied in computer vision.
Many challenges remain.
For instance, the classification of {critical configurations}, as carried out by~\cite{critical} in the setting of points viewed by projective cameras, gives precise conditions for when the 3D geometry cannot be uniquely recovered.
This has implications for the stability of reconstruction algorithms which are largely unexplored.
It is also desirable to better understand \emph{semialgebraic} constraints implicit in image formation, such as those which enforce that world points lie in front of the camera~\cite{chiral}. 
Two recent preprints that go in this direction are \cite{APST2020,PST2020,Max20}.
For other geometric problems in vision requiring robust estimation,
methods based on convex relaxations are popular alternatives to RANSAC, since solutions can often be certified by computing duality gaps.
The algebraic structures underlying these optimization problems are diverse, often drawing from the {sums of squares hierarchy} as in~\cite{kahl-henrion} or {synchronization} formulations as in~\cite{SESync,Arr16}.
We have not touched on many of the deeper tools from algebraic geometry appearing in algebraic vision, such as  in~\cite{Aho13,Aho14,Flo18,Lie20,Max20}. 
Nor have we addressed the role of nonlinear algebra in applications like photometric stereo~\cite{photo,Hea20}, methods for dynamic scenes~\cite{gpca}, or the wide variety of alternate camera models encountered in practice (e.g.~cameras with {radial distortion}~\cite{radial}).
These offer many potential directions for future research.
\section{Tensors and their decompositions}\label{section tensors}
\emph{by Paul Breiding}
\medskip

A \emph{tensor} $A$ is a $p$-dimensional array $A = (a_{i_1,\ldots,i_p})_{1\leq i_1\leq n_1,\ldots,1\leq i_p\leq n_p}$ filled with (real or complex) numbers $a_{i_1,\ldots,i_p}.$ The number $p$ is also called the \emph{order} of $A.$ For $p=2$ this gives the definition of a matrix. For $p\geq 3$ tensors are higher-dimensional analogues of matrices. In this section, we mainly consider real tensors.
The space of real $n_1\times \cdots \times n_p$-tensors is denoted by $\mathbb R^{n_1\times \cdots \times n_p}.$
Order-$3$ tensors can be visualized as cubes:
$$A \ =\  \mlatensor{2,3,2.5}  \,\in\,\mathbb R^{n_1\times n_2\times n_3}$$
While matrices are the central object in linear algebra, tensors define the field of {multilinear algebra}.

Of particular interest in applications are decompositions of a tensor $A,$ which reveal information about the data stored in $A$ \cite{PFS2016}. We seek to decompose $A$ as $A_1+\cdots+A_r,$ where the $A_i$ are some sort of \emph{simple} tensors. This means that the tensors~$A_i$ are elements in a low-dimensional model~$V,$ which usually is a real algebraic variety. In this section, we present two examples of this: the {Canonical Polyadic Decomposition} (CPD) and the {Block Term Decomposition} (BTD).

The CPD is the decomposition with simple tensors of the form $\mathbf{a}^{(1)} \otimes \cdots \otimes \mathbf{a}^{(p)} \coloneqq (a^{(1)}_{i_1}\cdots a^{(p)}_{i_p}),$ where the $\mathbf{a}^{(j)}=(a^{(j)}_i)\in\mathbb R^{n_j}$ are vectors. Tensors of the form $\mathbf{a}^{(1)} \otimes \cdots \otimes \mathbf{a}^{(p)}$ are called \emph{rank-one tensors}. The set of all rank-one tensors is a smooth projective algebraic variety, called \emph{Segre variety}.
\begin{definition}[CPD of tensors]\label{def_cpd}
Let $A\in\mathbb R^{n_1\times \cdots \times n_d}.$ We call a decomposition of the form
 \begin{equation}\label{cpd}
   A = \sum_{i=1}^r \mathbf{a}_i^{(1)} \otimes \cdots \otimes \mathbf{a}_i^{(p)},
 \end{equation}
where $\mathbf{a}_i^{(j)}\in\mathbb R^{n_j},$
a {\em canonical polyadic decomposition} (CPD) of $A.$ If $r$ is the minimal number such that we can write $A$ like in (\ref{cpd}), we say that $A$ has \emph{rank} $r.$
\end{definition}
\noindent For $p=3$ the CPD can be visualized as follows:
$$
\mlatensor{2,3,2.5} \ =\ \mlarankonetensor{2,3,2.5} \,+ \, \cdots \,+\, \mlarankonetensor{2,3,2.5}
$$
We can illustrate the meaning of this decomposition using an example from statistics: the relationship between random variables $X_1,\ldots,X_p,$ where $X_i$ can take $n_i$ states, can be recorded in a tensor $A=(a_{i_1,\ldots,i_p})\in \mathbb R^{n_1\times \cdots \times n_p}$ with $a_{i_1,\ldots,i_p}$ being the probability that $X_1=i_1, \ldots, X_p=i_p.$ The variables are independent if and only if the rank of $A$ is 1. The joint distribution is a mixture distribution of $r$ independent models if and only if the rank of $A$ is $r.$ Computing CPDs is also central in signal processing~\cite{7038247}, computational complexity
\cite{ottaviani2020tensor}, and other fields. The geometry of the CPD is well-studied in nonlinear algebra. One result worth mentioning is~\cite{DOT2018}, where the authors show cases in which the best rank-$r$ approximation of a generic tensor lies in the space spanned by its critical rank-$1$ approximations. For applications this implies that we can precondition the problem of computing best rank-$r$ approximations by first computing critical rank-$1$ approximations.

BTD is the decomposition, where the simple tensors are so-called \emph{low-multilinear rank} tensors. The definition of this is as follows: let $k=(k_1,\ldots,k_p)$ be a vector of integers with $1\leq k_i\leq n_i.$ We say that $A$ has multilinear-rank (at most)~$k,$ if there exists another tensor $S\in\mathbb R^{k_1\times \cdots  \times k_p},$ called the \emph{core tensor}, and matrices $U_i\in\mathbb R^{n_i\times k_i},$ called \emph{factor matrices}, such that we can write
$A = (U^{(1)}\otimes \cdots \otimes U^{(p)})\, S.$
Here, $U^{(1)}\otimes \cdots \otimes U^{(p)}$ is the linear map $\mathbb R^{k_1\times \cdots \times k_p} \to \mathbb R^{n_1\times \cdots \times n_p}$ defined by the action on rank-one tensors $(U^{(1)}\otimes \cdots \otimes U^{(p)})(\mathbf{x}^{(1)} \otimes \cdots \otimes \mathbf{x}^{(p)})\coloneqq (U^{(1)}\mathbf{x}^{(1)}) \otimes \cdots \otimes (U^{(p)}\mathbf{x}^{(p)}).$ Intuitively speaking, $(U^{(1)}\otimes \cdots \otimes U^{(p)})$ multiplies the $i$-th side of a tensor with $U^{(i)}$ as depicted in \eqref{btd_figure}. A~BTD decomposes $A$ as the sum of multilinear rank-$k$ tensors.
\begin{definition}[BTD of tensors]
Let $A\in\mathbb R^{n_1\times \cdots \times n_d}$ and $k=(k_1,\ldots,k_p).$ A decomposition
 \begin{equation}\label{btd}
   A = \sum_{i=1}^r (U^{(1)}_i\otimes \cdots \otimes U^{(p)}_i)\, S_i,
 \end{equation}
where $U_i^{(j)}\in \mathbb R^{n_j\times k_j}$ and $S_i\in \mathbb R^{k_1\times \cdots \times k_p},$ is called a
\emph{block term decomposition} (BTD) of $A.$ If $r$ is the minimal number such that we can write $A$ like in \eqref{btd}, we say that $A$ has BTD-\emph{rank} $r.$
\end{definition}
Note that CPD is the special case $k=(1,\ldots,1).$ For order-$3$ tensors a BTD can be visualized as follows:
\begin{equation}\label{btd_figure}\mlatensor{2,3,2.5} \ = \  \mlarankLMNtensor{2,3,2.5}{1,1,1}\,+\,\cdots \,+\, \mlarankLMNtensor{2,3,2.5}{1,1,1}\end{equation}
In the example from statistics, the CPD defined a tensor as a mixture of independence models. By contrast, BTD defines a tensor as a mixture of distributions which allow correlations between the variables. A prime example, where this is relevant, is detecting epileptic seizures
\cite{HCSPVHL2014}. The interaction between the variables in this case is extremely complex, so that a mixture of independence models is not the appropriate model.

For both CPD and BTD, and in general for any tensor decomposition, questions about uniqueness are important. These questions ask for formats and ranks, where a generic tensor $A$ of low rank has a \emph{unique} decomposition. In this case, we also say that $A$ is \emph{identifiable}. Note that for this we must have $p\geq 3.$ Matrices never have a unique CPD or BTD. These questions a priori fix the rank. The problem of computing or estimating ranks of tensors, also called \emph{model selection}
\cite{RKG2020}, is usually much more difficult.

Many tensor decompositions can be formulated within the following framework: let $V\subseteq \mathbb R^N$ be a real algebraic variety and $r\geq 2.$ Let $\phi_r: V^{\times r} \to \mathbb R^N, (\mathbf x_1,\ldots,\mathbf x_r)\mapsto \mathbf x_1+\cdots+\mathbf x_r.$
The image $\mathcal{J}_r(V)\coloneqq\phi_r(V^{\times r})$ is called the \emph{$r$-th join} of $V$. Similarly, we can define $\phi_r:
(V^{\mathbb C})^{\times r} \to \mathbb C^N$ and~$\mathcal{J}_r(V^{\mathbb C})$ for a complex algebraic variety $V^{\mathbb C}.$ For $\mathbf y\in\mathbb C^N$ we call the points in $\phi_r^{-1}(\mathbf y)$ \emph{rank-$r$ decompositions } of $y$. We say that $\mathbf y$ has a {\em unique} rank-$r$ decomposition, if $\phi_r^{-1}(\mathbf y)$ consists of only a single point $\mathbf x=(\mathbf x_1,\ldots,\mathbf x_r)$ modulo permutations of the $\mathbf x_i.$ We define the \emph{generic rank} $r_g(V^{\mathbb C})$ to be the smallest~$r,$ such that the {Zariski closure} of $\mathcal{J}_r(V^{\mathbb C})$ is $\mathbb C^N.$ The following lemma connects real and complex identifiability. Its assumptions are mild and hold for most tensor decompositions.
\begin{lemma}[{\cite[Lemma 28]{QPL2016}}]
Let $V = \{\mathbf x\in\mathbb R^N \mid f_1(\mathbf x)=\cdots=f_k(\mathbf x)=0\}$ be a real algebraic variety~that is also a cone; i.e., it is closed under scalar multiplication. We assume that $V$ is not contained in a hyperplane. Let~$V^{\mathbb C}= \{\mathbf x\in\mathbb C^N \mid f_1(\mathbf x)=\cdots=f_k(\mathbf x)=0\}$ be the corresponding complex variety and~$X\subseteq \mathbb P^{N-1}$ be the projective variety whose cone is $V^\mathbb{C}$; that is, $V^\mathbb{C} = \widehat{X}$ when using the notation from \Cref{section intro}. Let $r<r_g(V^{\mathbb C}).$ We assume that (a)~$X$ is irreducible, and (b)~that $V^{\mathbb C}$ has a real smooth point.
Then, if a generic point in~$\mathcal{J}_r(V^{\mathbb C})$ has a unique complex rank-r decomposition, a generic point in $\mathcal{J}_r(V)$ has a unique real rank-r decomposition.
\end{lemma}
The lemma shows that for studying identifiability of real tensors we can use methods from complex algebraic geometry. For CPD this was worked out in \cite[Section 5]{QPL2016} and \cite[Section 4]{BBV2019}.

Identifiability is relevant for applications, because it implies that the problem of computing tensor decompositions is {well-posed}: if we can compute one solution, we know that this solution is the right one. In the signal processing literature there exist several results on the uniqueness of the CPD, most of which rely on Kruskal's criterion~\cite{Kruskal1977}. This is a criterion for order-$3$ tensors, which relies on results from linear algebra. Nonlinear algebra goes beyond this. An overview about identifiability from the perspective of algebraic geometry is given by Angelini in~\cite{Angelini2017}. We also refer to the book of Landsberg~\cite{Landsberg2012}.

Uniqueness of BTDs is less well-studied than for CPDs. Results exist for the decomposition of tensors into tensors of multilinear rank $(1, k_1, k_2)$
\cite{Yang2014}, $(1, k, k)$
\cite{Lathauwer2011,DL2020}, $(k_1, k_2, 1)$
\cite{SL2015}, and $(k_1,k_2,k_3)$ \cite{Lathauwer2008DecompositionsOA}. Studying BTDs from the perspective of nonlinear algebra is a promising research direction, which would complement these results.
Another interesting research direction is studying the algebraic geometry of {structured tensor decompositions}. A prominent example of this is the symmetric CPD, i.e., CPDs of the form $A=\sum_{i=1}^r \mathbf a_i\otimes \mathbf a_i\otimes \cdots\otimes \mathbf a_i,$ where for each term all factors are equal. The identifiability of this decomposition is discussed in~\cite{COV2017} (surprisingly, it was recently shown by Shitov~\cite{Shitov2018} that the symmetric rank defined by this decomposition can differ from the rank in Definition \ref{def_cpd}).
Other examples of structured decompositions in the literature are {coupled CPDs}~\cite{SL2017, SL2017b,SL2015b} or
{Tensor Networks}~\cite{BDLG21}.

 {\footnotesize\linespread{0.8}
\bibliographystyle{abbrv}
\bibliography{bibliography}

\begin{thebibliography}{100}

\bibitem{AbeGri}
S.~Abenda and P.~G. Grinevich.
\newblock {Reducible {M}-curves for {L}e-networks in the totally-nonnegative
  {G}rassmannian and KP-II multiline soliton}.
\newblock {\em Selecta Math.}, 25:1--64, 2019.

\bibitem{Abl}
M.~J. Ablowitz.
\newblock Line soliton interactions.
\newblock Available at
  \url{https://sites.google.com/site/ablowitz/line-solitons}.

\bibitem{AblBal}
M.~J. Ablowitz and D.~E. Baldwin.
\newblock Nonlinear shallow ocean-wave soliton interactions on flat beaches.
\newblock {\em Phys. Rev. E (3)}, 86(036305), 2012.

\bibitem{ALSS}
M.~F. Adamer, A.~L. L\H{o}rincz, A.-L. Sattelberger, and B.~Sturmfels.
\newblock Algebraic {A}nalysis of {R}otation {D}ata.
\newblock {\em Algebraic Statistics}, 11(2):189--211, 2020.

\bibitem{AdlMos}
M.~Adler and J.~Moser.
\newblock On a class of polynomials connected with the {K}orteweg-de {V}ries
  equation.
\newblock {\em Comm. Math. Phys.}, 61:1--30, 1978.

\bibitem{Rome}
S.~Agarwal, Y.~Furukawa, N.~Snavely, I.~Simon, B.~Curless, S.~M. Seitz, and
  R.~Szeliski.
\newblock Building {R}ome in a day.
\newblock {\em Communications of the ACM}, 54(10):105--112, 2011.

\bibitem{APST2020}
S.~Agarwal, A.~Pryhuber, R.~Sinn, and R.~R. Thomas.
\newblock The {C}hiral {D}omain of a {C}amera {A}rrangement.
\newblock arXiv:2003.09265, 2020.

\bibitem{AgoAme}
D.~Agostini and C.~Am{\'e}ndola.
\newblock Discrete {G}aussian distributions via theta functions.
\newblock {\em SIAGA}, 2:1--30, 2019.

\bibitem{AgoCelEke}
D.~Agostini, T.~{\"{O}}. {\c{C}}elik, and D.~Eken.
\newblock Numerical reconstruction of curves from their {J}acobians.
\newblock arXiv:2103.03138, 2021.

\bibitem{AgoCelSruStu}
D.~Agostini, T.~{\"{O}}. {\c{C}}elik, J.~Struwe, and B.~Sturmfels.
\newblock Theta surfaces.
\newblock {\em Vietnam J. Math.}, 2020.

\bibitem{CelAgoStu}
D.~Agostini, T.~{\"{O}}. {\c{C}}elik, and B.~Sturmfels.
\newblock The {D}ubrovin threefold of an algebraic curve.
\newblock arXiv:2005.08244, 2020. To appear in {\em Nonlinearity}.

\bibitem{AgoChu}
D.~Agostini and L.~Chua.
\newblock Computing theta functions with {J}ulia.
\newblock {\em J. Softw. Algebra Geom.}, 11:41--51, 2021.

\bibitem{AgoFevManStu}
D.~Agostini, C.~Fevola, Y.~Mandelshtam, and B.~Sturmfels.
\newblock K{P} solitons from tropical limits.
\newblock arXiv:2101.10392, 2021.

\bibitem{Aho14}
C.~Aholt and L.~Oeding.
\newblock The ideal of the trifocal variety.
\newblock {\em Math. Comp.}, 83(289):2553--2574, 2014.

\bibitem{Aho13}
C.~Aholt, B.~Sturmfels, and R.~Thomas.
\newblock A {H}ilbert scheme in computer vision.
\newblock {\em Canad. J. Math.}, 65(5):961--988, 2013.

\bibitem{AEMS21}
R.~Ait El~Manssour and A.-L. Sattelberger.
\newblock Combinatorial {D}ifferential {A}lgebra of~$x^p$.
\newblock arXiv:2102.03182, 2021.

\bibitem{allman2010identifiability}
E.~S. Allman, S.~Petrovi{\'c}, J.~A. Rhodes, and S.~Sullivant.
\newblock Identifiability of two-tree mixtures for group-based models.
\newblock {\em IEEE/ACM transactions on computational biology and
  bioinformatics}, 8(3):710--722, 2010.

\bibitem{amendola2020computing}
C.~Am{\'e}ndola, L.~D. Garc{\'\i}a-Puente, R.~Homs, O.~Kuznetsova, and H.~J.
  Motwani.
\newblock Computing maximum likelihood estimates for {G}aussian graphical
  models with {M}acaulay2.
\newblock arXiv:2012.11572, 2020.

\bibitem{amendola2021invariant}
C.~Améndola, K.~Kohn, P.~Reichenbach, and A.~Seigal.
\newblock Invariant theory and scaling algorithms for maximum likelihood
  estimation.
\newblock {\em SIAM Journal on Applied Algebra and Geometry}, 5(2):304–337,
  2021.

\bibitem{Angelini2017}
E.~Angelini.
\newblock On complex and real identifiability of tensors.
\newblock {\em Rivista di Matematica della Universita di Parma}, 8:367--377, 12
  2017.

\bibitem{aoki2012markov}
S.~Aoki, H.~Hara, and A.~Takemura.
\newblock {\em Markov bases in algebraic statistics}, volume 199.
\newblock Springer Science \& Business Media, 2012.

\bibitem{Arr16}
F.~Arrigoni, B.~Rossi, and A.~Fusiello.
\newblock Spectral synchronization of multiple views in {S}{E} (3).
\newblock {\em SIAM Journal on Imaging Sciences}, 9(4):1963--1990, 2016.

\bibitem{Bertini}
D.~J. Bates, J.~D. Hauenstein, A.~J. Sommese, and C.~W. Wampler.
\newblock {Bertini: Software for Numerical Algebraic Geometry}.
\newblock Available at bertini.nd.edu with permanent doi:
  dx.doi.org/10.7274/R0H41PB5.

\bibitem{BelBobEnoItsMat}
E.~D. Belokolos, A.~I. Bobenko, V.~Z. Enolskii, A.~R. Its, and V.~B. Matveev.
\newblock {\em Algebro-geometric approach to non- linear integrable equations}.
\newblock Springer Series in Nonlinear Dynamics. Springer-Verlag, Berlin, 1994.

\bibitem{BBV2019}
C.~Beltrán, P.~Breiding, and N.~Vannieuwenhoven.
\newblock Pencil-based algorithms for tensor rank decomposition are not stable.
\newblock {\em SIAM J. Matrix Anal. Appl.}, 40, 02 2019.

\bibitem{BDLG21}
A.~Bernardi, C.~De~Lazzari, and F.~Gesmundo.
\newblock Dimension of tensor network varieties.
\newblock arXiv:2101.03148, 2021.

\bibitem{bernstein2017unimodular}
D.~I. Bernstein and S.~Sullivant.
\newblock Unimodular binary hierarchical models.
\newblock {\em J. Combin. Theory Ser. A}, 123:97--125, 2017.

\bibitem{bezanson2017julia}
J.~Bezanson, A.~Edelman, S.~Karpinski, and V.~B. Shah.
\newblock Julia: A fresh approach to numerical computing.
\newblock {\em SIAM review}, 59(1):65--98, 2017.

\bibitem{dickenstein:regions}
F.~Bihan, A.~Dickenstein, and M.~Giaroli.
\newblock Lower bounds for positive roots and regions of multistationarity in
  chemical reaction networks.
\newblock {\em J. Algebra}, 542:367--411, 2020.

\bibitem{Bj2}
J.-E. Bj{\"o}rk.
\newblock {\em Analytic {${\mathcal{D}}$}-{M}odules and {A}pplications}, volume
  247 of {\em Mathematics and its Applications}.
\newblock Kluwer Academic Publishers Group, Dordrecht, 1993.

\bibitem{Blekherman:Parrilo:Thomas:SDOptAndConvAlgGeo}
G.~Blekherman, P.~Parrilo, and R.~Thomas.
\newblock {\em Semidefinite optimization and convex algebraic geometry}.
\newblock Society for Industrial and Applied Mathematics, 2012.

\bibitem{boege2020reciprocal}
T.~Boege, J.~I. Coons, C.~Eur, A.~Maraj, and F.~R{\"o}ttger.
\newblock Reciprocal maximum likelihood degrees of {B}rownian motion tree
  models.
\newblock arXiv:2009.11849, 2020.

\bibitem{magma}
W.~Bosma, J.~Cannon, and C.~Playoust.
\newblock The {M}agma algebra system. {I}. {T}he user language.
\newblock {\em J. Symbolic Comput.}, 24(3-4):235--265, 1997.
\newblock Computational algebra and number theory (London, 1993).

\bibitem{HC.jl}
P.~Breiding and S.~Timme.
\newblock Homotopy{C}ontinuation.jl: {A} {P}ackage for {H}omotopy
  {C}ontinuation in {J}ulia.
\newblock In {\em Mathematical Software -- ICMS 2018}, pages 458--465, Cham,
  2018. Springer International Publishing.

\bibitem{BruSijZot}
N.~Bruin, J.~Sijsling, and A.~Zotine.
\newblock Numerical computation of endomorphism rings of {J}acobians.
\newblock In {\em 13th Algorithmic Number Theory Symposium}, volume~2 of {\em
  Open Book Series}, pages 155--171, 2019.

\bibitem{BMS}
C.~Bruschek, H.~Mourtada, and J.~Schepers.
\newblock Arc spaces and {R}ogers--{R}amanujan identities.
\newblock {\em Ramanujan J.}, 30:9--38, 2013.

\bibitem{BucEnoLey}
V.~Buchstaber, V.~Enolski, and D.~Leykin.
\newblock Rational analogs of abelian functions.
\newblock {\em Funct. Anal. Appl.}, 33:83--94, 1999.

\bibitem{Budur}
N.~Budur.
\newblock Bernstein--{S}ato ideals and local systems.
\newblock {\em Ann. Inst. Fourier (Grenoble)}, 65(2):549--603, 2015.

\bibitem{BMT}
N.~Budur, M.~Musta{\c t}{\u a}, and Z.~Teitler.
\newblock The monodromy conjecture for hyperplane arrangements.
\newblock {\em Geom. Dedicata}, 153(1):131--137, 2011.

\bibitem{bvwz}
N.~{Budur}, R.~{van der Veer}, L.~{Wu}, and P.~{Zhou}.
\newblock Zero loci of {B}ernstein--{S}ato ideals.
\newblock arXiv:1907.04010, 2019. To appear in {\em Invent. Math.}

\bibitem{caviness2012quantifier}
B.~F. Caviness and J.~R. Johnson.
\newblock {\em Quantifier elimination and cylindrical algebraic decomposition}.
\newblock Springer Science \& Business Media, 2012.

\bibitem{Cel}
T.~{\"{O}}. {\c{C}}elik.
\newblock {T}homae-{W}eber formula: Algebraic computations of theta constants.
\newblock {\em Int. Math. Res. Not. IMRN}, rnz302.

\bibitem{macis}
T.~{\"{O}}. {\c{C}}elik, A.~Jamneshan, G.~Mont\'{u}far, B.~Sturmfels, and
  L.~Venturello.
\newblock Optimal transport to a variety.
\newblock {\em Mathematical Aspects of Computer and Information Sciences,
  Springer Lecture Notes in Computer Science}, 11989:364--381, 2020.

\bibitem{CJMSV}
T.~{\"{O}}. {\c{C}}elik, A.~Jamneshan, G.~Montúfar, B.~Sturmfels, and
  L.~Venturello.
\newblock Wasserstein distance to independence models.
\newblock {\em J. Symbolic Comput.}, 104:855--873, 2021.

\bibitem{xiao:regularchain}
C.~Chen, J.~H. Davenport, M.~Moreno~Maza, B.~Xia, and R.~Xiao.
\newblock Computing with semi-algebraic sets represented by triangular
  decomposition.
\newblock In {\em Proceedings of the 2011 International Symposium on Symbolic
  and Algebraic Computation (ISSAC 2011)}, pages 75--82. ACM Press, 2011.

\bibitem{Hom4PSArticle}
T.~Chen, T.-L. Lee, and T.-Y. Li.
\newblock {Hom4PS-3: A Parallel Numerical Solver for Systems of Polynomial
  Equations Based on Polyhedral Homotopy Continuation Methods}.
\newblock In H.~Hong and C.~Yap, editors, {\em Mathematical Software -- ICMS
  2014}, pages 183--190. Springer Berlin Heidelberg, 2014.

\bibitem{COV2017}
L.~Chiantini, G.~Ottaviani, and N.~Vannieuwenhoven.
\newblock On generic identifiability of symmetric tensors of subgeneric rank.
\newblock {\em Trans. Amer. Math. Soc.}, 369:4021--4042, 06 2017.

\bibitem{CRS21}
Y.~Cid-Ruiz and B.~Sturmfels.
\newblock Primary {D}ecomposition with {D}ifferential {O}perators.
\newblock arXiv:2101.03643, 2021.

\bibitem{cifuentes2020sums}
D.~Cifuentes, T.~Kahle, and P.~Parrilo.
\newblock Sums of squares in {M}acaulay2.
\newblock {\em J. Softw. Algebra Geom}, 10(1):17--24, 2020.

\bibitem{collins1975quantifier}
G.~E. Collins.
\newblock Quantifier elimination for real closed fields by cylindrical
  algebraic decompostion.
\newblock In {\em Automata theory and formal languages}, pages 134--183.
  Springer, 1975.

\bibitem{Connelly-GenericGlobalRigidity}
R.~Connelly.
\newblock Generic {G}lobal {R}igidity.
\newblock {\em Discrete Comput. Geom.}, 33(4):549--563, 2005.

\bibitem{ConnellyServatiusHigher-orderRigidity}
R.~Connelly and H.~Servatius.
\newblock Higher-order rigidity---what is the proper definition?
\newblock {\em Discrete Comput. Geom.}, 11(2):193--200, 1994.

\bibitem{FeliuPlos}
C.~Conradi, E.~Feliu, M.~Mincheva, and C.~Wiuf.
\newblock Identifying parameter regions for multistationarity.
\newblock {\em PLoS Comput. Biol.}, 13(10):e1005751, 2017.

\bibitem{conradi-PNAS}
C.~Conradi, D.~Flockerzi, J.~Raisch, and J.~Stelling.
\newblock Subnetwork analysis reveals dynamic features of complex (bio)chemical
  networks.
\newblock {\em Proc. Nat. Acad. Sci. USA}, 104(49):19175--80, 2007.

\bibitem{shiu:hopf}
C.~Conradi, M.~Mincheva, and A.~Shiu.
\newblock Emergence of oscillations in a mixed-mechanism phosphorylation
  system.
\newblock {\em Bull. Math. Biol.}, 81(6):1829--1852, 2019.

\bibitem{coons2020generalized}
J.~I. Coons, J.~Cummings, B.~Hollering, and A.~Maraj.
\newblock Generalized cut polytopes for binary hierarchical models.
\newblock arXiv:2008.00043, 2020.

\bibitem{coons2021quasi}
J.~I. Coons and S.~Sullivant.
\newblock Quasi-independence models with rational maximum likelihood estimator.
\newblock {\em J. Symbolic Comput.}, 104:917--941, 2021.

\bibitem{coons2021toric}
J.~I. Coons and S.~Sullivant.
\newblock {Toric geometry of the {C}avender-{F}arris-{N}eyman model with a
  molecular clock}.
\newblock {\em Adv. in Appl. Math}, 123:102119, 2021.

\bibitem{Coutinho}
S.~C. Coutinho.
\newblock {\em A {P}rimer of {A}lgebraic $D$-{M}odules}, volume~33 of {\em
  Student Texts}.
\newblock London Mathematical Society, 1995.

\bibitem{CLO}
D.~Cox, J.~Little, and D.~O'Shea.
\newblock {\em Ideals, Varieties, and Algorithms: An Introduction to
  Computational Algebraic Geometry and Commutative Algebra}.
\newblock Undergraduate Texts in Mathematics. Springer New York, 4th edition,
  2007.

\bibitem{craciun2008}
G.~Craciun, J.~W. Helton, and R.~J. Williams.
\newblock Homotopy methods for counting reaction network equilibria.
\newblock {\em Math. Biosci}, 216(2):140--149, 2008.

\bibitem{Lathauwer2011}
L.~De~Lathauwer.
\newblock {Blind Separation of Exponential Polynomials and the Decomposition of
  a Tensor in Rank-$(L_r,L_r,1)$ Terms}.
\newblock {\em SIAM J. Matrix Anal. Appl.}, 32(4):1451--1474, 2011.

\bibitem{triang}
J.~A. De~Loera, J.~Rambau, and F.~Santos.
\newblock {\em Triangulations}, volume~25 of {\em Algorithms and computations
  in Mathematics}.
\newblock Springer-Verlag Berlin Heidelberg, 2010.

\bibitem{Singular}
W.~Decker, G.-M. Greuel, G.~Pfister, and H.~Sch\"onemann.
\newblock {\sc Singular} {4-1-3} --- {A} computer algebra system for polynomial
  computations.
\newblock \url{http://www.singular.uni-kl.de}, 2020.

\bibitem{algcurves}
B.~Deconinck, M.~Heil, A.~Bobenko, M.~van Hoeij, and M.~Schmies.
\newblock Computing {R}iemann theta functions.
\newblock {\em Math. Comp.}, 73:1417--1442, 2004.

\bibitem{Demazure}
M.~Demazure.
\newblock Sur deux problemes de reconstruction.
\newblock Technical Report RR-0882, INRIA, July 1988.
\newblock url: \url{https://hal.inria.fr/inria-00075672}.

\bibitem{derksen2020maximum}
H.~Derksen, V.~Makam, and M.~Walter.
\newblock Maximum likelihood estimation for tensor normal models via castling
  transforms.
\newblock arXiv:2011.03849, 2020.

\bibitem{diaconis1998algebraic}
P.~Diaconis, B.~Sturmfels, et~al.
\newblock Algebraic algorithms for sampling from conditional distributions.
\newblock {\em Ann. Statist.}, 26(1):363--397, 1998.

\bibitem{dickenstein2019multistationarity}
A.~Dickenstein, M.~P. Millan, A.~Shiu, and X.~Tang.
\newblock Multistationarity in structured reaction networks.
\newblock {\em Bull. Math. Biol.}, 81(5):1527--1581, 2019.

\bibitem{DM2007}
V.~Dolotin and A.~Morozov.
\newblock {\em {I}ntroduction to {N}on-linear {A}lgebra}.
\newblock World Scientific Publishing, Hackensack, NJ, 2007.

\bibitem{DL2020}
I.~Domanov and L.~D. Lathauwer.
\newblock {On Uniqueness and Computation of the Decomposition of a Tensor into
  Multilinear Rank-$(1,L_r,L_r)$ Terms}.
\newblock {\em SIAM J. Matrix Anal. Appl.}, 41(2):747--803, 2020.

\bibitem{control}
P.~Donnell, M.~Banaji, A.~Marginean, and C.~Pantea.
\newblock Control: an open source framework for the analysis of chemical
  reaction networks.
\newblock {\em Bioinformatics}, 30(11):1633--1634, 2014.

\bibitem{draisma2015noetherianity}
J.~Draisma, R.~Eggermont, R.~Krone, and A.~Leykin.
\newblock Noetherianity for infinite-dimensional toric varieties.
\newblock {\em Algebra Number Theory}, 9(8):1857--1880, 2015.

\bibitem{DHOST}
J.~Draisma, E.~Horobe\c{t}, G.~Ottaviani, B.~Sturmfels, and R.~R. Thomas.
\newblock The {E}uclidean distance degree of an algebraic variety.
\newblock {\em Found. Comput. Math.}, 16(1):99--149, 2016.

\bibitem{draisma2013groups}
J.~Draisma, S.~Kuhnt, P.~Zwiernik, et~al.
\newblock Groups acting on {G}aussian graphical models.
\newblock {\em Ann. Statist.}, 41(4):1944--1969, 2013.

\bibitem{DOT2018}
J.~Draisma, G.~Ottaviani, and A.~Tocino.
\newblock Best rank-k approximations for tensors: generalizing
  {E}ckart--{Y}oung.
\newblock {\em Res. Math. Sci.}, pages 5--27, 2018.

\bibitem{duarte2020equations}
E.~Duarte and C.~G{\"o}rgen.
\newblock Equations defining probability tree models.
\newblock {\em J. Symbolic Comput.}, 99:127--146, 2020.

\bibitem{DMS21}
E.~Duarte, O.~Marigliano, and B.~Sturmfels.
\newblock Discrete statistical models with rational maximum likelihood
  etimates.
\newblock {\em Bernoulli}, 27(1):135--154, 2021.

\bibitem{Dub}
B.~Dubrovin.
\newblock Theta functions and non-linear equations.
\newblock {\em Russian Math. Surveys}, 36:11--92, 1981.

\bibitem{PLMP}
T.~Duff, K.~Kohn, A.~Leykin, and T.~Pajdla.
\newblock {P}{L}{M}{P}-point-line minimal problems in complete multi-view
  visibility.
\newblock In {\em Proceedings of the IEEE International Conference on Computer
  Vision}, pages 1675--1684, 2019.

\bibitem{dye2020maximum}
S.~Dye, K.~Kohn, F.~Rydell, and R.~Sinn.
\newblock Maximum likelihood estimation for nets of conics.
\newblock arXiv:2011.08989, 2020.

\bibitem{Ehrenpreis}
L.~Ehrenpreis.
\newblock {\em Fourier Analysis in Several Complex Variables}, volume XVII of
  {\em Pure and Applied Mathematics}.
\newblock Wiley-Interscience Publishers, New York-London-Sydney, 1970.

\bibitem{EicZag}
M.~Eichler and D.~Zagier.
\newblock {\em The {T}heory of {J}acobi {F}orms}, volume~55 of {\em Progress in
  Mathematics}.
\newblock Birkh\"{a}user Boston, 1985.

\bibitem{Eisenbud}
D.~Eisenbud.
\newblock {\em Commutative Algebra: With a View Toward Algebraic Geometry}.
\newblock Graduate Texts in Mathematics. Springer New York, 2013.

\bibitem{crnttoolbox}
P.~Ellison, M.~Feinberg, H.~Ji, and D.~Knight.
\newblock Chemical reaction network toolbox, version 2.2.
\newblock Available online at \url{http://www.crnt.osu.edu/CRNTWin}, 2012.

\bibitem{TRPLP}
R.~Fabbri, T.~Duff, H.~Fan, M.~H. Regan, D.~d. C.~d. Pinho, E.~Tsigaridas,
  C.~W. Wampler, J.~D. Hauenstein, P.~J. Giblin, B.~Kimia, A.~Leykin, and
  T.~Pajdla.
\newblock {T}{R}{P}{L}{P} - {T}rifocal {R}elative {P}ose from {L}ines at
  {P}oints.
\newblock In {\em Proceedings of the IEEE/CVF Conference on Computer Vision and
  Pattern Recognition (CVPR)}, June 2020.

\bibitem{Feinbergss}
M.~Feinberg.
\newblock The existence and uniqueness of steady states for a class of chemical
  reaction networks.
\newblock {\em Arch. Rational Mech. Anal.}, 132(4):311--370, 1995.

\bibitem{feinberg_2019}
M.~Feinberg.
\newblock {\em Foundations of Chemical Reaction Network Theory}, volume 202 of
  {\em Applied Mathematical Sciences}.
\newblock Springer International Publishing, 2019.

\bibitem{feliu_newinj}
E.~Feliu.
\newblock Injectivity, multiple zeros, and multistationarity in reaction
  networks.
\newblock {\em Proc. Roy. Soc. Edinburgh Sect. A}, doi:10.1098/rspa.2014.0530,
  2014.

\bibitem{MyCRN1}
E.~Feliu, N.~Kaihnsa, T.~de~Wolff, and O.~Y{\"u}r{\"u}k.
\newblock The kinetic space of multistationarity in dual phosphorylation.
\newblock {\em J. Dynam. Differential Equations}, 2020.

\bibitem{fevola2020pencils}
C.~Fevola, Y.~Mandelshtam, and B.~Sturmfels.
\newblock {Pencils of quadrics: Old and new}.
\newblock arXiv:2009.04334, 2020.

\bibitem{RANSAC}
M.~A. Fischler and R.~C. Bolles.
\newblock Random sample consensus: a paradigm for model fitting with
  applications to image analysis and automated cartography.
\newblock {\em Communications of the ACM}, 24(6):381--395, 1981.

\bibitem{Flo18}
G.~Fl\o~ystad, J.~Kileel, and G.~Ottaviani.
\newblock The {C}how form of the essential variety in computer vision.
\newblock {\em J. Symbolic Comput.}, 86:97--119, 2018.

\bibitem{matlabtheta}
J.~Frauendiener, C.~Jaber, and C.~Klein.
\newblock Efficient computation of multidimensional theta functions.
\newblock {\em J. Geom. Phys.}, 127:147--158, 2019.

\bibitem{garcia2005algebraic}
L.~D. Garcia, M.~Stillman, and B.~Sturmfels.
\newblock Algebraic geometry of {B}ayesian networks.
\newblock {\em J. Symbolic Comput.}, 39(3-4):331--355, 2005.

\bibitem{Gau}
P.~Gaudry.
\newblock Fast genus 2 arithmetic based on {T}heta functions.
\newblock {\em J. Math. Cryptol. J}, 3(1):243--265, 2007.

\bibitem{geiger2006toric}
D.~Geiger, C.~Meek, B.~Sturmfels, et~al.
\newblock On the toric algebra of graphical models.
\newblock {\em Ann. Statist.}, 34(3):1463--1492, 2006.

\bibitem{GKZ}
I.~M. Gelfand, M.~M. Kapranov, and A.~V. Zelevinsky.
\newblock {\em Discriminants, resultants, and multidimensional determinants}.
\newblock Mathematics: Theory \& Applications. Birkh\"{a}user Boston, Inc.,
  Boston, MA, 1994.

\bibitem{gorgen2021staged}
C.~G{\"{o}}rgen, A.~Maraj, and L.~Nicklasson.
\newblock Staged tree models with toric structure.
\newblock 2107.04516.

\bibitem{GLS}
P.~G\"orlach, C.~Lehn, and A.-L. Sattelberger.
\newblock {A}lgebraic {A}nalysis of the {H}ypergeometric {F}unction ${_1F_1}$
  of a {M}atrix {A}rgument.
\newblock {\em Beitr\"age Algebra Geom.}, 62:397--427, 2021.

\bibitem{GortlerHealyThurston}
S.~J. Gortler, A.~D. Healy, and D.~P. Thurston.
\newblock Characterizing generic global rigidity.
\newblock {\em Amer. J. Math.}, 132(4):897--939, 2010.

\bibitem{vBR}
H.-C. Graf~von Bothmer and K.~Ranestad.
\newblock {A general formula for the algebraic degree in semidefinite
  programming}.
\newblock {\em Bull. Math. Biol.}, 41(2):193--197, 02 2009.

\bibitem{GraverServatius-CombinatorialRigidity}
J.~Graver, B.~Servatius, and H.~Servatius.
\newblock {\em Combinatorial rigidity}, volume~2 of {\em Graduate Studies in
  Mathematics}.
\newblock American Mathematical Society, Providence, RI, 1993.

\bibitem{M2}
D.~R. Grayson and M.~E. Stillman.
\newblock Macaulay2, a software system for research in algebraic geometry.
\newblock Available at \url{http://www.math.uiuc.edu/Macaulay2/}, 2020.

\bibitem{SingularCommAlg}
G.~Greuel, G.~Pfister, O.~Bachmann, C.~Lossen, and H.~Sch{\"o}nemann.
\newblock {\em A Singular Introduction to Commutative Algebra}.
\newblock Springer Nature Book Archives Millennium. Springer, 2002.

\bibitem{pollefeys-multicamera}
C.~H{\"a}ne, L.~Heng, G.~H. Lee, F.~Fraundorfer, P.~Furgale, T.~Sattler, and
  M.~Pollefeys.
\newblock 3{D} visual perception for self-driving cars using a multi-camera
  system: {C}alibration, mapping, localization, and obstacle detection.
\newblock {\em Image and Vision Computing}, 68:14--27, 2017.

\bibitem{harkonen2020algebraic}
M.~H{\"a}rk{\"o}nen, B.~Hollering, F.~T. Kashani, and J.~I. Rodriguez.
\newblock Algebraic optimization degree.
\newblock {\em ACM Commun. Comput. Algebra}, 54(2):44--48, 2020.

\bibitem{critical}
R.~Hartley and F.~Kahl.
\newblock Critical configurations for projective reconstruction from multiple
  views.
\newblock {\em Int. J. Comput. Vis.}, 71(1):5--47, 2007.

\bibitem{HZ-2003}
R.~Hartley and A.~Zisserman.
\newblock {\em Multiple View Geometry in Computer Vision}.
\newblock Cambridge, 2nd edition, 2003.

\bibitem{chiral}
R.~I. Hartley.
\newblock Chirality.
\newblock {\em Int. J. Comput. Vis.}, 26(1):41--61, 1998.

\bibitem{Hartshorne}
R.~Hartshorne.
\newblock {\em Algebraic Geometry}.
\newblock Graduate Texts in Mathematics. Springer New York, 2013.

\bibitem{HNTT13}
H.~Hashiguchi, Y.~Numata, N.~Takayama, and A.~Takemura.
\newblock The holonomic gradient method for the distribution function of the
  largest root of a {W}ishart matrix.
\newblock {\em J. Multivariate Anal.}, 117:296--312, 2013.

\bibitem{HL2015}
J.~D. Hauenstein and A.~C. Liddell.
\newblock A hybrid symbolic-numeric approach to exceptional sets of generically
  zero-dimensional systems.
\newblock In {\em Proceedings of the 2015 International Workshop on Parallel
  Symbolic Computation}, PASCO '15, page 53–60. Association for Computing
  Machinery, 2015.

\bibitem{Hea20}
K.~Heal, J.~Wang, S.~J. Gortler, and T.~Zickler.
\newblock A {L}ighting-{I}nvariant {P}oint {P}rocessor for {S}hading.
\newblock pages 94--102, 2020.

\bibitem{catastrophe}
A.~Heaton and S.~Timme.
\newblock Catastrophe in {E}lastic {T}ensegrity {F}rameworks.
\newblock arXiv:2009.13408, 2020.

\bibitem{Hilbert:1888}
D.~Hilbert.
\newblock {\"{U}ber die Darstellung definiter Formen als Summe von
  Formenquadraten}.
\newblock {\em Math. Ann.}, 32(3):342--350, 1888.

\bibitem{hillar2012finite}
C.~J. Hillar and S.~Sullivant.
\newblock Finite {G}r{\"o}bner bases in infinite dimensional polynomial rings
  and applications.
\newblock {\em Adv. Math.}, 229(1):1--25, 2012.

\bibitem{Hir}
R.~Hirota.
\newblock {\em The Direct Method in Soliton Theory}.
\newblock Cambridge Tracts in Mathematics. Cambridge University Press, 2004.

\bibitem{HitSegWar}
N.~Hitchin, G.~Segal, and R.~S. Ward.
\newblock {\em Integrable Systems: Twistors, Loop Groups, and Riemann
  Surfaces}.
\newblock Oxford Graduate Texts in Mathematics 4. Clarendon Press, Oxford,
  1999.

\bibitem{hollering2021identifiability}
B.~Hollering and S.~Sullivant.
\newblock Identifiability in phylogenetics using algebraic matroids.
\newblock {\em J. Symbolic Comput.}, 104:142--158, 2021.

\bibitem{hocsten2002grobner}
S.~Ho{\c{s}}ten and S.~Sullivant.
\newblock Gr{\"o}bner bases and polyhedral geometry of reducible and cyclic
  models.
\newblock {\em J. Combin. Theory Ser. A}, 100(2):277--301, 2002.

\bibitem{HTT08}
R.~Hotta, K.~Takeuchi, and T.~Tanisaki.
\newblock {\em {$D$}-{M}odules, {P}erverse {S}heaves, and {R}epresentation
  {T}heory}, volume 236 of {\em Progress in Mathematics}.
\newblock Birkh\"{a}user Boston, Inc., Boston, MA, 2008.
\newblock Translated from the 1995 Japanese edition by Takeuchi.

\bibitem{huh2013varieties}
J.~Huh.
\newblock Varieties with maximum likelihood degree one.
\newblock {\em J. Algebr. Stat.}, 5:1--17, 2014.

\bibitem{HS14}
J.~Huh and B.~Sturmfels.
\newblock {Likelihood Geometry}.
\newblock In {\em Combinatorial algebraic geometry}, volume 2108 of {\em
  Lecture notes in mathematics}, pages 63--117. Springer, New York, 2014.

\bibitem{HCSPVHL2014}
B.~Hunyadi, D.~Camps, L.~Sorber, W.~V. Paesschen, M.~D. Vos, S.~V. Huffel, and
  L.~D. Lathauwer.
\newblock Block term decomposition for modelling epileptic seizures.
\newblock {\em EURASIP J. Adv. Signal Process.}, (1):139, 2014.

\bibitem{Hurwitz:AMGM}
A.~Hurwitz.
\newblock Ueber den {V}ergleich des arithmetischen und des geometrischen
  {M}ittels.
\newblock {\em J. Reine Angew. Math.}, 108:266--268, 1891.

\bibitem{Iliman:deWolff:Circuits}
S.~Iliman and T.~de~Wolff.
\newblock Amoebas, nonnegative polynomials and sums of squares supported on
  circuits.
\newblock {\em Res. Math. Sci.}, 3(9), 2016.

\bibitem{joshi2015survey}
B.~Joshi and A.~Shiu.
\newblock A survey of methods for deciding whether a reaction network is
  multistationary.
\newblock {\em Math. Model. Nat. Phenom.}, 10(5):47--67, 2015.

\bibitem{kahl-henrion}
F.~Kahl and D.~Henrion.
\newblock Globally optimal estimates for geometric reconstruction problems.
\newblock {\em Int. J. Comput. Vis.}, 74(1):3--15, 2007.

\bibitem{Kashiwara}
M.~Kashiwara.
\newblock {$B$}-functions and holonomic systems. {R}ationality of roots of
  {$B$}-functions.
\newblock {\em Invent. Math.}, 38:33--53, 1976.

\bibitem{Kas84}
M.~Kashiwara.
\newblock The {R}iemann--{H}ilbert problem for holonomic systems.
\newblock {\em Publ. Res. Inst. Math. Sci.}, 20(2):319--365, 1984.

\bibitem{orealgebra}
M.~Kauers, M.~Jaroschek, and F.~Johansson.
\newblock Ore polynomials in {S}age.
\newblock ar{X}iv:1306.4263, 2013.

\bibitem{Kil17}
J.~Kileel.
\newblock Minimal problems for the calibrated trifocal variety.
\newblock {\em SIAM J. Appl. Algebra Geom.}, 1(1):575--598, 2017.

\bibitem{Kod}
Y.~Kodama.
\newblock {\em K{P} Solitons and the {G}rassmannians. Combinatorics and
  Geometry and Two-dimensional Wave Patterns}, volume~22 of {\em Briefs in
  Mathematical Physics}.
\newblock Springer, 2017.

\bibitem{KodWil}
Y.~Kodama and L.~Williams.
\newblock {KP solitons and total positivity for the Grassmannian}.
\newblock {\em Invent. Math.}, 198:637--699, 2014.

\bibitem{KodXie}
Y.~Kodama and Y.~Xie.
\newblock {Space curves and solitons of the KP hierarchy: I. The l-th
  generalized {K}d{V} hierarchy}.
\newblock {\em SIGMA}, 17(024), 2021.

\bibitem{Kolchin}
E.~R. Kolchin.
\newblock {\em Differential {A}lgebra and {A}lgebraic {G}roups}.
\newblock Pure and applied mathematics. Academic Press, 1973.

\bibitem{CK}
C.~Koutschan.
\newblock {HolonomicFunctions} (user's guide).
\newblock Technical Report 10-01, RISC Report Series, Johannes Kep\-ler
  University, Linz, Austria, 2010.

\bibitem{Koyama}
T.~Koyama.
\newblock The annihilating ideal of the {F}isher integral.
\newblock {\em Kyushu J. Math}, 74:415--427, 2020.

\bibitem{Kri1977}
I.~M. Krichever.
\newblock Methods of algebraic geometry in the theory of nonlinear equations.
\newblock {\em Russian Math. Surveys}, 32(6):44--48, 1977.

\bibitem{krone2017hilbert}
R.~Krone, A.~Leykin, and A.~Snowden.
\newblock Hilbert series of symmetric ideals in infinite polynomial rings via
  formal languages.
\newblock {\em J. Algebra}, 485:353--362, 2017.

\bibitem{Kruskal1977}
J.~B. Kruskal.
\newblock Three-way arrays: rank and uniqueness of trilinear decompositions,
  with application to arithmetic complexity and statistics.
\newblock {\em Linear Algebra Appl.}, 18:95--138, 1977.

\bibitem{AutoMin}
Z.~Kukelova, M.~Bujnak, and T.~Pajdla.
\newblock Automatic generator of minimal problem solvers.
\newblock In {\em European Conference on Computer Vision}, pages 302--315.
  Springer, 2008.

\bibitem{radial}
Z.~Kukelova and T.~Pajdla.
\newblock A minimal solution to the autocalibration of radial distortion.
\newblock In {\em 2007 IEEE Conference on Computer Vision and Pattern
  Recognition}, pages 1--7. IEEE, 2007.

\bibitem{periods}
P.~Lairez.
\newblock Computing periods of rational integrals.
\newblock {\em Math. Comp.}, 85:1719--1752, 2016.

\bibitem{LMS19}
P.~Lairez, M.~Mezzarobba, and M.~Safey El~Din.
\newblock {Computing the volume of compact semi-algebraic sets}.
\newblock In {\em {ISSAC 2019 - International Symposium on Symbolic and
  Algebraic Computation}}, Beijing, China, July 2019. {ACM}.

\bibitem{Laman1970}
G.~Laman.
\newblock On graphs and rigidity of plane skeletal structures.
\newblock {\em J. Engrg. Math.}, 4:331--340, 1970.

\bibitem{Landsberg2012}
J.~M. Landsberg.
\newblock {\em Tensors: Geometry and applications}.
\newblock Graduate Studies in Mathematics. AMS, Providence, Rhode Island, 2012.

\bibitem{syz}
V.~Larsson, K.~Astrom, and M.~Oskarsson.
\newblock Efficient solvers for minimal problems by syzygy-based reduction.
\newblock In {\em Proceedings of the IEEE Conference on Computer Vision and
  Pattern Recognition}, pages 820--829, 2017.

\bibitem{Laserre:MomentsPosPoly}
J.~B. Lasserre.
\newblock {\em Moments, positive polynomials and their applications}, volume~1
  of {\em Imperial College Press Optimization}.
\newblock World Scientific, 2010.

\bibitem{Lathauwer2008DecompositionsOA}
L.~D. Lathauwer.
\newblock Decompositions of a higher-order tensor in block terms - {Part II}:
  Definitions and uniqueness.
\newblock {\em SIAM J. Matrix Anal. Appl.}, 30:1033--1066, 2008.

\bibitem{laurent1999}
M.~Laurent and N.~Kellershohn.
\newblock Multistability: a major means of differentiation and evolution in
  biological systems.
\newblock {\em Trends Biochem. Sciences}, 24(11):418--422, 1999.

\bibitem{lauritzen2019maximum}
S.~Lauritzen, C.~Uhler, P.~Zwiernik, et~al.
\newblock Maximum likelihood estimation in {G}aussian models under total
  positivity.
\newblock {\em Ann. Statist.}, 47(4):1835--1863, 2019.

\bibitem{dmodlib}
V.~Levandovskyy and J.~Mart\'{i}n-Morales.
\newblock dmod\_lib: A {{\tt Singular:Plural}} library for algorithms for
  algebraic {$D$}-modules.
\newblock \url{https://www.singular.uni-kl.de/Manual/4-2-0/sing\_537}.

\bibitem{NumericalAlgebraicGeometryArticle}
A.~Leykin.
\newblock Numerical {A}lgebraic {G}eometry for {M}acaulay2.
\newblock {\em J. Softw. Algebra Geom.}, 3:5--10, 2011.

\bibitem{Lie20}
M.~Lieblich and L.~Van~Meter.
\newblock Two {H}ilbert schemes in computer vision.
\newblock {\em SIAM J. Appl. Algebra Geom.}, 4(2):297--321, 2020.

\bibitem{Lit}
J.~Little.
\newblock Translation manifolds and the converse to {A}bel’s theorem.
\newblock {\em Compos. Math}, 49:147--171, 1983.

\bibitem{Liu}
Q.~Liu.
\newblock {\em Algebraic Geometry and Arithmetic Curves}.
\newblock Oxford graduate texts in mathematics. Oxford University Press, 2006.

\bibitem{OnTropComp}
M.~Luxton and Z.~Qu.
\newblock Some results on tropical compactifications.
\newblock {\em Trans. Amer. Math. Soc.}, 363(9):4853--4876, 2011.

\bibitem{MS15}
D.~Maclagan and B.~Sturmfels.
\newblock {\em {Introduction to Tropical Geometry}}, volume 161 of {\em
  Graduate studies in mathematics}.
\newblock American Mathematical Society, Providence, R.I., 2015.

\bibitem{Mai16}
P.~Maisonobe.
\newblock Filtration relative, l'id\'eal de {B}ernstein et ses pentes.
\newblock hal-01285562v2, 2016.

\bibitem{Maple}
{Maplesoft, a division of Waterloo Maple Inc.}
\newblock Maple, 2019.

\bibitem{maraj2020algebraic}
A.~Maraj.
\newblock {\em Algebraic and Geometric Properties of Hierarchical Models}.
\newblock PhD thesis, 2020.

\bibitem{maraj2019equivariant}
A.~Maraj and U.~Nagel.
\newblock Equivariant {H}ilbert series for hierarchical models.
\newblock {\em Algebraic Statistics}, 12(1):21--42, 2021.

\bibitem{Marshall:Book}
M.~Marshall.
\newblock {\em Positive polynomials and sums of squares}, volume 146 of {\em
  Mathematical Surveys and Monographs}.
\newblock American Mathematical Society, 2008.

\bibitem{mathrepo}
{Max Planck Institute for Mathematics in the Sciences}.
\newblock Mathrepo. {M}athematical data and software.
\newblock Repository website of the MPI MiS,
  \url{https://mathrepo.mis.mpg.de/}.

\bibitem{Max20}
L.~G. Maxim, J.~I. Rodriguez, and B.~Wang.
\newblock Euclidean distance degree of the multiview variety.
\newblock {\em SIAM J. Appl. Algebra Geom.}, 4(1):28--48, 2020.

\bibitem{Meb84}
Z.~Mebkhout.
\newblock Une \'{e}quivalence de cat\'{e}gories.
\newblock {\em Compos. Math.}, 51:51--62, 1984.

\bibitem{MS2018}
M.~Micha\l{}ek and B.~Sturmfels.
\newblock Lecture: {I}ntroduction to {N}onlinear {A}lgebra.
\newblock Available at \url{https://youtu.be/1EryuvBLY80}, 2018.

\bibitem{MS2020}
M.~Micha\l{}ek and B.~Sturmfels.
\newblock {\em {I}nvitation to {N}onlinear {A}lgebra}, volume 211 of {\em
  Graduate Studies in Mathematics}.
\newblock American Mathematical Society, 2021.
\newblock Available for download at
  https://personal-homepages.mis.mpg.de/michalek/NonLinearAlgebra.pdf.

\bibitem{combalg}
E.~Miller and B.~Sturmfels.
\newblock {\em Combinatorial commutative algebra}, volume 227 of {\em Graduate
  Texts in Mathematics}.
\newblock Springer-Verlag, New York, 2005.

\bibitem{misra2020gaussian}
P.~Misra and S.~Sullivant.
\newblock {G}aussian graphical models with toric vanishing ideals.
\newblock {\em Ann. Inst. Statist. Math.}, pages 1--29, 2020.

\bibitem{Mui70}
R.~J. Muirhead.
\newblock Systems of partial differential equations for hypergeometric
  functions of matrix argument.
\newblock {\em Ann. Math. Statist.}, 41:991--1001, 1970.

\bibitem{Mui82}
R.~J. Muirhead.
\newblock {\em Aspects of multivariate statistical theory}.
\newblock John Wiley \& Sons, Inc., New York, 1982.
\newblock Wiley Series in Probability and Mathematical Statistics.

\bibitem{Mum1}
D.~Mumford.
\newblock {\em {T}ata {L}ectures on {T}heta {I}}.
\newblock Modern Birkh{\"a}user Classics. Birkh\"{a}user Boston, 2007.

\bibitem{nagel2017equivariant}
U.~Nagel and T.~R{\"o}mer.
\newblock Equivariant {H}ilbert series in non-{N}oetherian polynomial rings.
\newblock {\em J. Algebra}, 486:204--245, 2017.

\bibitem{nagpal2021symmetric}
R.~Nagpal and A.~Snowden.
\newblock Symmetric ideals of the infinite polynomial ring.
\newblock arXiv:2107.13027, 2021.

\bibitem{NNNOSTT11}
H.~Nakayama, K.~Nishiyama, M.~Noro, K.~Ohara, T.~Sei, N.~Takayama, and
  A.~Takemura.
\newblock Holonomic gradient descent and its application to the
  {F}isher--{B}ingham integral.
\newblock {\em Adv. in Appl. Math.}, 47(3):639--658, 2011.

\bibitem{Nak2010}
A.~Nakayashiki.
\newblock On {A}lgebraic {E}xpansions of {S}igma {F}unctions for $(n, s)$
  curves.
\newblock {\em Asian J. Math.}, 14:175--212, 2010.

\bibitem{Naka}
A.~Nakayashiki.
\newblock Degeneration of trigonal curves and solutions of the {KP}-hierarchy.
\newblock {\em Nonlinearity}, 31:3567--3590, 2018.

\bibitem{NR}
J.~Nie and K.~Ranestad.
\newblock Algebraic degree of polynomial optimization.
\newblock {\em SIAM J. Optim.}, 20(1):485--502, 2009.

\bibitem{NRS}
J.~Nie, K.~Ranestad, and B.~Sturmfels.
\newblock The algebraic degree of semidefinite programming.
\newblock {\em Math. Program.}, 122(2, Ser. A):379--405, 2010.

\bibitem{Nister}
D.~Nist{\'e}r.
\newblock An efficient solution to the five-point relative pose problem.
\newblock {\em IEEE transactions on pattern analysis and machine intelligence},
  26(6):756--770, 2004.

\bibitem{Noro}
M.~Noro.
\newblock System of partial differential equations for the hypergeometric
  function ${_1F_1}$ of a matrix argument on diagonal regions.
\newblock {\em {ISSAC} '16: {P}roceedings of the {ACM} on {I}nternational
  {S}ymposium of {S}ymbolic and {A}lgebraic {C}omputation}, pages 381--388,
  July 2016.

\bibitem{NovManPitaZak}
S.~Novikov, S.~Manakov, L.~Pitaevskii, and V.~Zakharov.
\newblock {\em Theory of Solitons: The Inverse Scattering Method}.
\newblock Monographs in Contemporary Mathematics. Springer, 1984.

\bibitem{Ollivier}
F.~Ollivier.
\newblock Standard bases of differential ideals.
\newblock {\em Lecture Notes in Comput. Sci.}, 508:304--321, 1990.

\bibitem{ottaviani2020tensor}
G.~Ottaviani and P.~Reichenbach.
\newblock Tensor rank and complexity.
\newblock arXiv:2004.01492, 2020.

\bibitem{OtSo}
G.~Ottaviani and L.~Sodomaco.
\newblock The distance function from a real algebraic variety.
\newblock {\em Comput. Aided Geom. Design}, 82:101927, 2020.

\bibitem{ozbudak2004}
E.~M. Ozbudak, M.~Thattai, H.~N. Lim, B.~I. Shraiman, and A.~Van~Oudenaarden.
\newblock Multistability in the lactose utilization network of escherichia
  coli.
\newblock {\em Nature}, 427(6976):737--740, 2004.

\bibitem{Palamodov}
V.~P. Palamodov.
\newblock {\em Linear Differential Qperators with Constant Coefficients},
  volume 168 of {\em Grundlehren der mathematischen Wissenschaften}.
\newblock Springer-Verlag, New York-Berlin, 1970.

\bibitem{PFS2016}
E.~E. Papalexakis, C.~Faloutsos, and N.~D. Sidiropoulos.
\newblock Tensors for data mining and data fusion: Models, applications, and
  scalable algorithms.
\newblock {\em ACM Trans. Intell. Syst. Technol.}, 8(2), Oct. 2016.

\bibitem{PerezMillan}
M.~P{\'e}rez~Mill{\'a}n, A.~Dickenstein, A.~Shiu, and C.~Conradi.
\newblock Chemical reaction systems with toric steady states.
\newblock {\em Bull. Math. Biol.}, 74:1027--1065, 2012.

\bibitem{petrovic2019markov}
S.~Petrovi{\'c}.
\newblock What is... a markov basis?
\newblock {\em Notices of the American Mathematical Society}, 66(7), 2019.

\bibitem{Pistone1996GeneralisedCW}
G.~Pistone and H.~Wynn.
\newblock Generalised confounding with grobner bases.
\newblock {\em Biometrika}, 83:653--666, 1996.

\bibitem{Pollaczek-Geiringer1927}
H.~Pollaczek-Geiringer.
\newblock \"{U}ber die {G}liederung ebener {F}achwerke.
\newblock {\em J. Appl. Math. Mech.}, 7:58--72, 1927.

\bibitem{PST2020}
A.~Pryhuber, R.~Sinn, and R.~R. Thomas.
\newblock Existence of two view chiral reconstructions.
\newblock arXiv:2011.07197, 2020.

\bibitem{QPL2016}
Y.~Qi, P.~Comon, and L.-H. Lim.
\newblock Semialgebraic geometry of nonnegative tensor rank.
\newblock {\em SIAM J. Matrix Anal. Appl.}, 37:1556--1580, 11 2016.

\bibitem{R}
{R Core Team}.
\newblock {\em {\tt R}: {A} Language and Environment for Statistical
  Computing}.
\newblock R Foundation for Statistical Computing, Vienna, Austria, 2013.

\bibitem{RegSte}
O.~Regev and N.~Stephens-Davidowitz.
\newblock An inequality for {G}aussians on lattices.
\newblock {\em SIAM J. Discrete Math}, 31(2):749--757, 2017.

\bibitem{surveyGKZ}
T.~Reichelt, M.~Schulze, C.~Sevenheck, and U.~Walther.
\newblock Algebraic aspects of hypergeometric differential equations.
\newblock {\em Beitr\"age Algebra Geom.}, February 2021.

\bibitem{Ritt}
J.~F. Ritt.
\newblock {\em Differential Algebra}, volume~14 of {\em American {M}athematical
  {S}ociety: {C}olloquium publications}.
\newblock American {M}athematical {S}ociety, 1950.

\bibitem{RKG2020}
A.~A. Rontogiannis, E.~Kofidis, and P.~Giampouras.
\newblock Block-term tensor decomposition: Model selection and computation.
\newblock arXiv:2002.09759, 2020.

\bibitem{SESync}
D.~M. Rosen, L.~Carlone, A.~S. Bandeira, and J.~J. Leonard.
\newblock {S}{E}-{S}ync: A certifiably correct algorithm for synchronization
  over the special {E}uclidean group.
\newblock {\em The International Journal of Robotics Research},
  38(2-3):95--125, 2019.

\bibitem{SabbahBS}
C.~Sabbah.
\newblock Proximit\'{e} \'{e}vanescente. {I}. {L}a structure polaire d’un
  {$D$}-module.
\newblock {\em Compositio Math.}, 62(3):283--328, 1987.

\bibitem{SST00}
M.~Saito, B.~Sturmfels, and N.~Takayama.
\newblock {\em Gr\"{o}bner deformations of hypergeometric differential
  equations}, volume~6 of {\em Algorithms and Computation in Mathematics}.
\newblock Springer-Verlag, Berlin, 2000.

\bibitem{gfun}
B.~Salvy and P.~Zimmermann.
\newblock {{\tt gfun}: A Maple Package for the Manipulation of Generating and
  Holonomic Functions in One Variable}.
\newblock {\em ACM Trans. Math. Software}, 20(2):163--177, 1994.

\bibitem{Sato}
M.~Sato.
\newblock Soliton equations as dynamical systems on infinite dimensional
  {G}rassmann manifold.
\newblock In {\em Nonlinear Partial Differential Equations in Applied Science;
  Proceedings of The U.S.-Japan Seminar, Tokyo, 1982}, North-Holland
  Mathematics Studies, pages 259--271. North-Holland, 1983.

\bibitem{SatStu19}
A.-L. Sattelberger and B.~Sturmfels.
\newblock {{$D$}-Modules and Holonomic Functions}.
\newblock arXiv:1910.01395, 2019.

\bibitem{TropicalMLE}
A.-L. Sattelberger and R.~van~der Veer.
\newblock {Maximum Likelihood Estimation from a Tropical and a
  {B}ernstein--{S}ato Perspective}.
\newblock arXiv:2101.03570, 2021.

\bibitem{SegWil}
G.~Segal and G.~Wilson.
\newblock Loop groups and equations of {K}d{V} type.
\newblock {\em Publ. Math. Inst. Hautes Études Sci.}, 61:5--65, 1985.

\bibitem{FarGruMan}
H.~Segal, S.~Grushevsky, and R.~S. Manni.
\newblock An explicit solution to the weak {S}chottky problem.
\newblock {\em Algebr. Geom.}, 8(3):358--373, 2021.

\bibitem{SSTOT13}
T.~Sei, H.~Shibata, A.~Takemura, K.~Ohara, and N.~Takayama.
\newblock Properties and applications of {F}isher distribution on the rotation
  group.
\newblock {\em J. Multivariate Anal.}, 116:440--455, 2013.

\bibitem{Shitov2018}
Y.~Shitov.
\newblock A counterexample to {C}omon's conjecture.
\newblock {\em SIAM Journal on Applied Algebra and Geometry}, 2(3):428--443,
  2018.

\bibitem{7038247}
N.~Sidiropoulos, L.~De~Lathauwer, X.~Fu, K.~Huang, E.~Papalexakis, and
  C.~Faloutsos.
\newblock Tensor decomposition for signal processing and machine learning.
\newblock {\em IEEE Transactions on Signal Processing}, 65(13):3551--3582,
  2017.

\bibitem{Sil}
R.~Silhol.
\newblock The {S}chottky problem for real genus 3 {M}-curves.
\newblock {\em Math. Z.}, 236:841--881, 2001.

\bibitem{SitharamStJohnSidman-HandbookGeometricConstraintSystems}
M.~Sitharam, A.~St.~John, and J.~Sidman, editors.
\newblock {\em Handbook of geometric constraint systems principles}.
\newblock Discrete Mathematics and its Applications (Boca Raton). CRC Press,
  Boca Raton, FL, 2019.

\bibitem{Sommese:Wampler:2005}
A.~J. Sommese and C.~W. Wampler.
\newblock {\em The {N}umerical {S}olution of {S}ystems of {P}olynomials
  {A}rising in {E}ngineering and {S}cience}.
\newblock World Scientific, 2005.

\bibitem{SL2015}
M.~Sorensen and L.~De~Lathauwer.
\newblock {Coupled Canonical Polyadic Decompositions and (Coupled)
  Decompositions in Multilinear Rank-$(L_r,n,L_r,n,1)$ Terms---Part I:
  Uniqueness}.
\newblock {\em SIAM J. Matrix Anal. Appl.}, 36(2):496--522, 2015.

\bibitem{SL2017}
M.~{Sorensen} and L.~{De Lathauwer}.
\newblock {Multidimensional Harmonic Retrieval via Coupled Canonical Polyadic
  Decomposition—Part I: Model and Identifiability}.
\newblock {\em IEEE Trans. Signal Process.}, 65(2):517--527, 2017.

\bibitem{SL2017b}
M.~{Sorensen} and L.~{De Lathauwer}.
\newblock {Multidimensional Harmonic Retrieval via Coupled Canonical Polyadic
  Decomposition—Part II: Algorithm and Multirate Sampling}.
\newblock {\em IEEE Trans. Signal Process.}, 65(2):528--539, 2017.

\bibitem{SL2015b}
D.~L.~L. Sorensen~M., Domanov~I.
\newblock {Coupled Canonical Polyadic Decompositions and (Coupled)
  Decompositions in Multilinear rank-$(L_r,n,L_r,n,1)$ terms --- Part II:
  Algorithms}.
\newblock {\em SIAM J. Matrix Anal. Appl.}, 36(3):1015--1045, 2015.

\bibitem{Strang-computationalscienceandengineering}
G.~Strang.
\newblock {\em Computational science and engineering}.
\newblock Wellesley-Cambridge Press, Wellesley, MA, 2007.

\bibitem{sturmfels1996grobner}
B.~Sturmfels.
\newblock {\em {Gr\"obner Bases and Convex Polytopes}}, volume~8 of {\em
  University Lecture Series}.
\newblock American Mathematical Soc., 1996.

\bibitem{Sturmfels02solvingsystems}
B.~Sturmfels.
\newblock {\em {Solving Systems of Polynomial Equations}}.
\newblock Number~97 in CBMS Regional Conferences Series. American Mathematical
  Society, 2002.

\bibitem{whatisGrobner}
B.~Sturmfels.
\newblock What is \ldots a {G}r\"obner basis?
\newblock {\em Notices Amer. Math. Soc.}, 52(10), 2005.

\bibitem{sturmfels2005toric}
B.~Sturmfels and S.~Sullivant.
\newblock Toric ideals of phylogenetic invariants.
\newblock {\em J. Comput. Biol.}, 12(2):204--228, 2005.

\bibitem{ST20}
B.~Sturmfels and S.~Telen.
\newblock {Likelihood Equations and Scattering Amplitudes}.
\newblock arXiv:2012.05041, 2020.

\bibitem{SU2010}
B.~Sturmfels and C.~Uhler.
\newblock Multivariate {G}aussian, semidefinite matrix completion, and convex
  algebraic geometry.
\newblock {\em Ann. Inst. Statist. Math.}, 62(4):603--638, 2010.

\bibitem{Sturmfels2020BrownianMT}
B.~Sturmfels, C.~Uhler, and P.~Zwiernik.
\newblock Brownian motion tree models are toric.
\newblock {\em Kybernetika}, 56:1154--1175, 2020.

\bibitem{sullivant2009gaussian}
S.~Sullivant.
\newblock Gaussian conditional independence relations have no finite complete
  characterization.
\newblock {\em J. Pure Appl. Algebra}, 213(8):1502--1506, 2009.

\bibitem{algstat2018}
S.~Sullivant.
\newblock {\em Algebraic statistics}, volume 194 of {\em Graduate Studies in
  Mathematics}.
\newblock American Mathematical Society, Providence, RI, 2018.

\bibitem{abelfunctions}
C.~Swierczewski and B.~Deconinck.
\newblock Computing {R}iemann theta functions in {S}age with applications.
\newblock {\em Math. Comput. Simulation}, 127:263--272, 2016.

\bibitem{Szeliski}
R.~Szeliski.
\newblock {\em Computer Vision - Algorithms and Applications}.
\newblock Texts in Computer Science. Springer, 2011.

\bibitem{hgmR}
N.~Takayama, T.~Koyama, T.~Sei, H.~Nakayama, and K.~Nishiyama.
\newblock {\em {\tt hgm}: Holonomic Gradient Method and Gradient Descent},
  2017.
\newblock {\tt R} package version 1.17.

\bibitem{sagemath}
{The Sage Developers}.
\newblock {\em {S}ageMath, the {S}age {M}athematics {S}oftware {S}ystem}.
\newblock \url{https://www.sagemath.org}.

\bibitem{torres:stability}
A.~Torres and E.~Feliu.
\newblock Detecting parameter regions for bistability in reaction networks.
\newblock arXiv:1909.13608, 2019.

\bibitem{Vakil}
R.~Vakil.
\newblock The {R}ising {S}ea: {F}oundations of {A}lgebraic {G}eometry.
\newblock Available at
  \url{http://math.stanford.edu/~vakil/216blog/index.html}.

\bibitem{PHCpack}
J.~Verschelde.
\newblock Algorithm 795: {PHC}pack: {A} {G}eneral-{P}urpose {S}olver for
  {P}olynomial {S}ystems by {H}omotopy {C}ontinuation.
\newblock {\em ACM Trans. Math. Softw.}, 25(2):251–276, June 1999.

\bibitem{gpca}
R.~Vidal, Y.~Ma, and S.~Sastry.
\newblock Generalized principal component analysis ({G}{P}{C}{A}).
\newblock {\em IEEE transactions on pattern analysis and machine intelligence},
  27(12):1945--1959, 2005.

\bibitem{wiuf-feliu}
C.~Wiuf and E.~Feliu.
\newblock Power-law kinetics and determinant criteria for the preclusion of
  multistationarity in networks of interacting species.
\newblock {\em SIAM J. Appl. Dyn. Syst.}, 12:1685--1721, 2013.

\bibitem{photo}
R.~J. Woodham.
\newblock Photometric method for determining surface orientation from multiple
  images.
\newblock In {\em Shape from shading}, pages 513--531. 1989.

\bibitem{Xiong:2003jt}
W.~Xiong and J.~E. Ferrell~Jr.
\newblock {A positive-feedback-based bistable 'memory module' that governs a
  cell fate decision}.
\newblock {\em Nature}, 426(6965):460--465, 2003.

\bibitem{Yang2014}
M.~Yang.
\newblock On partial and generic uniqueness of block term tensor
  decompositions.
\newblock {\em Annali dell'universita' di Ferrara}, 60(2):465--493, 2014.

\bibitem{Zei90}
D.~Zeilberger.
\newblock A holonomic systems approach to special functions identities.
\newblock {\em J. Comput. Appl. Math.}, 32(3):321--368, 1990.

\bibitem{Zobnin}
A.~I. Zobnin.
\newblock One-element differential standard bases with respect to inverse
  lexicographical orderings.
\newblock {\em J.~Math. Sci. (N. Y.)}, 163(5):523--533, 2009.

\end{thebibliography}
}

\bigskip

\bigskip

\noindent \textbf{Acknowledgements:}

\smallskip \small

\noindent Paul Breiding is funded by the Deutsche Forschungsgemeinschaft (DFG, German Research Foundation) -- Projektnummer 445466444. T\"urk\"u \"Ozl\"{u}m \c{C}elik is supported by Turkish Scientific and Technological Research Council (TÜB\.{I}TAK) -- TÜB\.{I}TAK 2236, Project number 1119B362000396. Timothy Duff is supported by a NSF Mathematical Sciences Postdoctoral Research Fellowship (DMS-2103310.)  Alex Heaton is supported by the Fields Institute for Research in Mathematical Sciences. Lorenzo Venturello is supported by the G\"oran Gustafsson foundation.

\bigskip

\bigskip

\linespread{1.2}
\noindent \textbf{Authors' addresses:}

\smallskip \small

\noindent Paul Breiding, MPI MiS, Leipzig\hfill paul.breiding@mis.mpg.de\\
T\"urk\"u \"Ozl\"um \c{C}elik, Bo\u{g}azi\c{c}i University, \.{I}stanbul\hfill  turkuozlum@gmail.com\\
Timothy F. Duff, University of Washington, Seattle \hfill timduff@uw.edu\\
Alex Heaton, Lawrence University, Appleton \hfill alexheaton2@gmail.com\\
Aida Maraj, University of Michigan, Ann Arbor \hfill  maraja@umich.edu\\
Anna-Laura Sattelberger, MPI MiS, Leipzig \hfill anna-laura.sattelberger@mis.mpg.de\\
Lorenzo Venturello, KTH Royal Institute of Technology, Stockholm \hfill lven@kth.se\\
O\u{g}uzhan Y\"ur\"uk, TU Berlin, Chair of Discrete Mathematics/Geometry
\hfill yuruk@math.tu-berlin.de
\end{document}